\documentclass[oneside, 12pt]{amsart}
\usepackage{amscd, amssymb, amsmath, mathrsfs}
\usepackage[english]{babel}
\usepackage{url}
\usepackage{tikz-cd}
\usepackage{booktabs}
\usepackage{dynkin-diagrams}
\usepackage[pdftex, colorlinks=true,  citecolor=blue, linkcolor=blue, linktocpage=true]{hyperref}

\DeclareUrlCommand\arXiv{\urlstyle{same}}

\newcommand\mylabel[1]{\label{#1}\marginpar{\vspace{-1ex}\medskip\medskip\footnotesize \tt #1}}
\renewcommand\mylabel[1]{\label{#1}}
\newcommand{\mydate}{
\number\day\space
\ifcase\month \or January\or February\or March\or April\or May\or June\or July\or August\or September\or October\or November\or December\fi 
\space\number\year}

\newtheorem{theorem}{Theorem}[section]
\newtheorem{maintheorem}{Theorem}
\newtheorem{lemma}[theorem]{Lemma}
\newtheorem{proposition}[theorem]{Proposition}
\newtheorem{corollary}[theorem]{Corollary}

\theoremstyle{definition}

\newtheorem{remark}[theorem]{Remark}

\newtheorem*{acknowledgement}{Acknowledgement}

\theoremstyle{remark}

\newcommand{\FF}{\mathbb{F}}

\newcommand{\PP}{\mathbb{P}}

\newcommand{\GG}{\mathbb{G}}

\newcommand{\ideala}{\mathfrak{a}}
\newcommand{\idealb}{\mathfrak{b}}

\newcommand{\idealh}{\mathfrak{h}}

\newcommand{\shM}{\mathscr{M}}

\newcommand{\shN}{\mathscr{N}}
\newcommand{\shL}{\mathscr{L}}

\newcommand{\Ad}{\operatorname{Ad}}
\newcommand{\ad}{\operatorname{ad}}

\newcommand{\Ass}{\operatorname{Ass}}
\newcommand{\Aut}{\operatorname{Aut}}

\newcommand{\Bl}{\operatorname{Bl}}

\newcommand{\Cokernel}{\operatorname{Coker}}

\newcommand{\depth}{\operatorname{depth}}

\newcommand{\Ext}{\operatorname{Ext}}

\newcommand{\Frac}{\operatorname{Frac}}

\newcommand{\GL}{\operatorname{GL}}
\newcommand{\liegl}{\operatorname{\mathfrak{gl}}}

\newcommand{\Grass}{\operatorname{Grass}}

\newcommand{\Hom}{\operatorname{Hom}}

\newcommand{\id}{{\operatorname{id}}}
\newcommand{\Image}{\operatorname{Im}}

\newcommand{\Kernel}{\operatorname{Ker}}

\newcommand{\length}{\operatorname{length}}
\newcommand{\Lie}{\operatorname{Lie}}

\newcommand{\lra}{\longrightarrow}

\newcommand{\Mat}{\operatorname{Mat}}

\newcommand{\maxid}{\mathfrak{m}}

\renewcommand{\O}{\mathscr{O}}

\newcommand{\pd}{\operatorname{pd}}

\newcommand{\PGL}{\operatorname{PGL}}

\newcommand{\quadand}{\quad\text{and}\quad}

\newcommand{\ra}{\rightarrow}

\newcommand{\red}{{\operatorname{red}}}
\newcommand{\Reg}{\operatorname{Reg}}

\newcommand{\Sing}{\operatorname{Sing}}
\newcommand{\SL}{\operatorname{SL}}
\newcommand{\liesl}{\operatorname{\mathfrak{sl}}}

\newcommand{\Spec}{\operatorname{Spec}}

\newcommand{\uHom}{\underline{\operatorname{Hom}}}

\newcommand{\lieg}{\mathfrak{g}}

\newcommand{\lieh}{\mathfrak{h}}

\newcommand{\tS}{{\tilde{S}}}

\newcommand{\tY}{{\tilde{Y}}}
\newcommand{\tf}{{\tilde{f}}}
\newcommand{\inert}{\text{\rm inert}}
\newcommand{\liepgl}{\operatorname{\mathfrak{pgl}}}

\begin{document}

\title[Surfaces of general type]
      {Surfaces of general type and $\mathfrak{sl}_2$-triples}

\author[Stefan Schr\"oer]{Stefan Schr\"oer}
\address{Heinrich Heine University D\"usseldorf, Faculty of Mathematics and Natural Sciences, Mathematical Institute, 40204 D\"usseldorf, Germany}
\curraddr{}
\email{schroeer@math.uni-duesseldorf.de}

\author[]{Nikolaos Tziolas}
\address{Department of Mathematics and Statistics,
University of Cyprus,
P.O.\ Box 20537, Nicosia, Cyprus}
\curraddr{}
\email{tziolas@ucy.ac.cy}



\begin{abstract}
The  $\liesl_2$-triples play a fundamental role for the structure theory of Lie algebras, and representation theory in general.
Here we investigate $\liesl_2$-triples of global vector fields on schemes $X$ in positive characteristics $p>0$,
and develop a general theory for actions of the corresponding height-one group scheme $G=\SL_2[F]$.
Sending a point to the Lie algebra of its stabilizer defines rational maps to various Grassmann varieties.
For surfaces of general type, this yields fibrations in  curves of genus $g\geq 2$ over the projective line.
Using properties  of the corresponding moduli stack $\shM_g$, we prove that there are no smooth surfaces of general type
with an $\liesl_2$-triple.
On the other hand, employing   Lefschetz pencils and Frobenius pullbacks  we show that canonical surfaces of general type
with such triples exist in abundance. In this connection, we classify the rational double points where the tangent sheaf is free
or the evaluation pairing with K\"ahler differentials is surjetive, including characteristic two. 
\end{abstract}

\maketitle
\tableofcontents

\section*{Introduction}
\mylabel{Introduction}

The  main goal of this paper is to understand infinitesimal symmetries of surfaces of general type in positive characteristics.
The general context is as follows:

Let $X$ be a proper scheme over a ground field 
$k$, for the moment  of arbitrary characteristic $p\geq 0$. Then $G=\Aut_{X/k}$ exists as a group scheme, 
the Lie algebra is the space $\lieg=H^0(X,\Theta_{X/k})$ of global vector fields, and
the connected component $G^0$ of the origin is of finite type.
Brion and the first author \cite{Brion; Schroeer 2023} showed that   every connected group scheme of finite type   arises
as some $\Aut^0_{X/k}$.
On the other hand, for smooth varieties $X$ of general type,   $\Aut_{X/k}$ is actually finite,    according to 
a   result of Martin-Deschamps and Lewin-M\'en\'egaux
\cite{Martin-Deschamps; Lewin-Menegaux 1978}.
Thus $\Aut^0_{X/k}$ is a singleton, and the story in characteristic zero ends here.
For $p>0$, however,  this finite group scheme might be non-reduced, but not  much seems to be known
about its structure.
In light of this state of affairs it is natural to concentrate attention on  the 
Frobenius kernel $G[F]=\Aut_{X/k}^0[F]$, or equivalently on the restricted Lie algebra $\lieg=H^0(X,\Theta_{X/k})$,
and try to understand that.

From now on  we are in characteristic  $p>0$. 
The above problem is vacuous in dimension $n=1$, in other words for curves   of genus at least two, because then the tangent sheaf 
has strictly negative degree. The situation changes drastically in dimension $n=2$:
The first surfaces of general type with non-zero global vector fields where constructed by Raynaud \cite{Raynaud 1978}.
Further results in this direction where obtained by Shepherd-Barron \cite{Shepherd-Barron 1991}  and the second author
\cite{Tziolas 2017}, \cite{Tziolas 2022}.
We showed in \cite{Schroeer; Tziolas 2023} that for smooth surfaces   of general type, and their canonical models   alike, 
among the startling multitude of restricted Lie algebras 
only $k^n$ and $ k^n\rtimes \liegl_1(k)$ and $ \liesl_2(k)$
are possible; the corresponding height-one group schemes are the Frobenius kernels
of the vector group $\GG_a^n$, the semidirect product  $\GG_a^n\rtimes\GG_m$, and the special linear group 
$\SL_2$, respectively. For all of them we constructed example,  with the nerve-wrecking exception of $\liesl_2(k)$.

The restricted Lie algebra $\liesl_2(k)$ is three-dimensional,  and plays a fundamental role
in the structure theory of Lie algebras and representation theory in general (confer \cite{LIE 7-8}, Chapter VIII, \S 1 and \S 11). 
However, applications in algebraic geometry are surprisingly scarce. The
traceless matrices $e=(\begin{smallmatrix} 0&1\\0&0\end{smallmatrix})$ and $h=(\begin{smallmatrix} 1&0\\0&-1\end{smallmatrix})$ and 
$f=(\begin{smallmatrix} 0&0\\-1&0\end{smallmatrix})$ form a basis, with  bracket and $p$-map   given by 
$$
[h,e]=2e,\quad [h,f]=-2f,\quad [e,f]=-h\quadand e^{[p]}=f^{[p]}=0,\quad h^{[p]}=h.
$$ 
For any restricted Lie algebra $\lieh$,  the injective homomorphisms $\liesl_2(k)\ra \lieh$   correspond to linearly independent
triples $(e,h,f)$ from   $\lieh$ satisfying the above relations. Following Dynkin and Bourbaki  we call them
\emph{$\liesl_2$-triples}.
We also find it convenient to say that a scheme $X$ \emph{admits an $\liesl_2$-triple} if  the restricted Lie algebra $\lieh=H^0(X,\Theta_{X/k})$ of global vector fields
contains such an $(e,h,f)$. Geometrically speaking, this  is nothing but a faithful action of the height-one group scheme $\SL_2[F]$. Note that by~\cite{Schroeer; Tziolas 2023}, a variety of foliation rank one has at most one, up to automorphisms of $\mathfrak{sl}_2(k)$, faithful $\mathfrak{sl}_2$-triple.

Our  first  main result is:

\begin{maintheorem}
(See Thm.\ \ref{thm:general type no triple}) No proper smooth surface  $X$ of general type admits an $\liesl_2$-triple.
\end{maintheorem}

To establish this, we  develop a general theory for   actions of the height-one group scheme $G=\SL_2[F]$
on integral schemes $Y$ whose generic stabilizer has order  $|G_\eta|=p^2$. Writing $\lieg=\liesl_2(k)$, we then consider  the so-called \emph{inertia map}
$$
f_\inert:Y\dashrightarrow \Grass^2(\lieg),\quad y\longmapsto \Lie(G_y)
$$
to the Grassmann variety of two-dimensional subspaces. This of course factors over the scheme of two-dimensional subalgebras,
which turn out to be the curve $B\subset \Grass^2(\lieg)$ defined by a quadratic equation $T_0^2-4T_1T_2=0$.
From this one can already deduces that the characteristic must be $p\neq 2$.  Now the key insight is:

\begin{maintheorem}
(See Thm.\ \ref{thm:resolution of fixed points}) 
On some $G$-modification $\tilde{Y}$, the action is fixed point free, and the foliation   
$\lieg\cdot\O_\tY\subset\Theta_{\tY/k}$ is a direct summand  
identified with $\tf^*_\inert(\Theta_{B/k})$.
\end{maintheorem}

For smooth surfaces $X$ of general type this leads, at least on some modification, to a pencil of smooth curves of genus $g\geq 2$.
By Szpiro's generalization \cite{Szpiro 1979} of  results of Parshin \cite{Parshin 1968} and Arakelov \cite{Arakelov 1971},
the classifying morphism $\PP^1\ra\shM_g$ to the moduli stack must be constant.
But this implies $X=C\times\PP^1$, which is impossible for surfaces of general type.

However, the story does not end here: 
Recall that each \emph{smooth surface    of general type} $S$ has   a unique    model $Y$ 
where the dualizing sheaf $\omega_Y$ is ample and the singularities $\O_{Y,y}$ are rational double points.
Let us call such $Y$  \emph{canonical surfaces of general type}.  
From a certain perspective they are of superior interest,  because they  form an Artin stack with finite inertia groups.
It is also convenient to have a name for 
their modifications $X\ra Y$ whose singularities $\O_{X,x}$ stay  rational double points; we propose to call them
\emph{RDP surfaces of general type}.
Note that Blanchard's Lemma gives an inclusion $\Aut^0_{X/k}\subset\Aut^0_{Y/k}$, which may or may not be
an equality. Our third main results shows  that such surfaces with $\liesl_2$-triples are surprisingly common:

\begin{maintheorem}
(See Thm.\ \ref{thm:existence with higher height}) 
Let $S$ be a smooth surface of general type in characteristic $p\geq 3$.
Then there is a purely inseparable alteration $X\ra S$ by some RDP surface of general type
having an $\liesl_2$-triple. Moreover, one may assume that the tangent sheaf $\Theta_{X/k}$ is locally free,
and  the action of $\SL_2[F]$ has no fixed points.
\end{maintheorem}

The idea is to use  \emph{Lefschetz pencils}  $C_t\subset S$, which are parameterized by $t\in\PP^1$ and 
becomes a fibration upon  blowing-up the axis $Z$; one then obtains the desired RDP surface of general type 
$X=\Bl_Z(S)\times_{\PP^1}(\PP^1,F)$ as  a \emph{relative Frobenius base-change of the   Lefschetz fibration}.
The $\SL_2[F]$-action arises from the canonical action on the second factor. 
Note that using higher Frobenius powers $F^n$ produces actions of iterated Frobenius kernels. Also note
that the  fibration $X\ra\PP^1$ can be seen as  morphism $\PP^1\ra\bar{\shM}_g$ into the Delinge--Mumford stack
of stable curves which factors over the relative Frobenius of the relative Frobenius map of the projective line. 
We believe that the  phenomenon of  inseparability in classifying maps to moduli stacks deserve further study.

The above construction produces  rational double points whose tangent modules $\Theta_{R/k}$
are free and   the evaluation pairing $\Theta_{R/k}\otimes\Omega^1_{R/k}\ra R$ is surjective.
In light of Artin's classification of rational double points \cite{Artin 1977}, this raises the question 
which rational double points have one or both of these properties.
Such question go back to Lipman  \cite{Lipman 1965}, and arise in many different contexts
(\cite{Wahl 1975},  \cite{Ekedahl; Hyland; Shepherd-Barron 2012}, 
\cite{Schroeer 2008}, \cite{Hirokado 2019}, \cite{Schroeer 2021}, \cite{Graf 2022}, \cite{Kawakami 2022}, \cite{Liedtke 2024}, \cite{Martin; Stadlmayr 2024}).
Combining the theory of \emph{minimal free resolutions} with techniques from \emph{Gr\"obner bases}, we indeed settle this,
including the most challenging case of characteristic two:

\begin{maintheorem}
(See Section \ref{sec:rational double points}) 
A rational double point in characteristic $p>0$ has free tangent module $\Theta_{R/k}$ 
or surjective evaluation pairing  $\Theta_{R/k}\otimes\Omega^1_{R/k}\ra R$
if and only if occurs in  table \ref{tab:behaviour rdp}.
\end{maintheorem}

\begin{table}
$$
\begin{array}{*{4}{l}}
\toprule
\text{\rm RDP}		& \text{\rm condition}	& \text{\rm tangent module free} 	& \text{\rm pairing surjective}\\
\toprule
A_l			& \text{\rm $l\equiv -1$ mod $p$}	& \text{\rm yes}			& \text{\rm yes}\\
\midrule
D_{2n}^0,D_{2n+1}^0	& p=2			& \text{\rm yes}			& \text{\rm yes}\\
D_{2n}^1,\ldots,D_{2n}^{n-1}&  			& \text{\rm yes}			& \text{\rm no}\\
D_{2n+1}^1,\ldots,D_{2n+1}^{n-1}& 		& \text{\rm no}				& \text{\rm no}\\
\midrule
E_8^0			& p=2,3,5		& \text{\rm yes}			& \text{\rm yes}\\
E_6^0,E_7^0		& p=2,3			& \text{\rm yes}			& \text{\rm yes}\\
E_6^1,E_7^1,E_8^1,E_8^2	& p=2			& \text{\rm yes}			& \text{\rm no}\\
\bottomrule
\end{array}
$$
\medskip
\caption{Properties of rational double points in positive characteristic}
\label{tab:behaviour rdp}
\end{table} 

Note that for RDP surfaces  of general type with faithful $\liesl_2$-triple, only $p\geq 3$ matters.
It would be interesting to see whether   $E_8^0\; (p=3,5)$ and $E_6^0,E_7^0\; (p=3)$ really occur  in
constructions like $X=\Bl_Z(S)\times_{\PP^1}(\PP^1,F)$, perhaps with certain pencils that violate the Lefschetz condition.


\emph{The paper is organized as follows:} In Section \ref{sec:generalities} we discuss some relevant facts on   automorphism group schemes,
surfaces of general type, rational double points, and $\liesl_2(k)$ as restricted Lie algebra.
In Section \ref{sec:pencils and triples} we analyze fibrations $f:X\ra \PP^1$ that are equivariant 
with respect to actions of the height-one group scheme
$G=\SL_2[F]$. Our core observations appear in Section \ref{sec:inertia map}, where we analyze the stabilizers for $G$-actions on schemes $X$,
and the resulting rational maps to Grassmann varieties of subalgebras  in $\lieg=\liesl_2(k)$.
These results are applied to smooth surfaces of general type in Section \ref{sec:surfaces}, where the above rational maps are
interpreted in term of the moduli stack $\shM_g$.  
Section \ref{sec:lefschetz pencils} contains construction of  canonical surfaces of general type having $\liesl_2$-triples.
In  Section \ref{sec:rational double points} we use homological algebra and Gr\"obner bases techniques to 
determine which rational double points have free tangent module
$\Theta_{R/k}$ or where the evaluation pairing   $\Theta_{R/k}\otimes\Omega^1_{R/k}\ra R$ is surjective.
In  the final    Section \ref{sec:rational double points with triples} we classify the rational double points 
admitting an  $\liesl_2$-triples and study when they are  fixed point free.

\begin{acknowledgement}
The research was supported by the Deutsche Forschungsgemeinschaft via the  project \emph{Varieties with Free Tangent Sheaves}, project
number 536205323. It  was  also conducted       in the framework of the   research training group
\emph{GRK 2240: Algebro-Geometric Methods in Algebra, Arithmetic and Topology}.
The first author also likes to thank the Department of Mathematics and Statistics  at the University of Cyprus for its hospitality during
a visit in October 2025. We also wish to thank Yuya Matsumoto for valuable comments. We also thank the referee for many useful comments,
which helped to improve the paper. In particular for the suggestion to
include ruled surfaces in Section 4, and to consider Lefschetz pencils
on general smooth surfaces in Section 5
\end{acknowledgement}

\section{Generalities}
\mylabel{sec:generalities}

Let $k$ be a ground field, for the moment of arbitrary characteristic $p\geq 0$,
and $X$ be a proper scheme. Matsumura and Oort \cite{Matsumura; Oort 1967} established that $\Aut_{X/k}$ is representable by 
a group scheme whose  connected component $\Aut^0_{X/k}$ of the origin is of finite type.
The Lie algebra is the space $H^0(X,\Theta_{X/k})$ of global vector fields, where  
$\Theta_{X/k}=\uHom(\Omega^1_{X/k},\O_X)$ denotes the \emph{tangent sheaf}. This is a  coherent sheaf,
which is locally free provided that $X$ is smooth, and than   rank equals   dimension.
Note that in positive characteristics there are also  singular schemes having locally free tangent sheaves.

Let us  say that a morphism $f:X\ra Y$  between proper schemes is   \emph{in Stein factorization}
if the canonical map $\O_Y\ra f_*(\O_X)$ is an isomorphism. Then Blanchard's Lemma (\cite{Blanchard 1956}, Proposition 1.1, 
see \cite{Brion 2017}, Theorem 7.2.1 for the   scheme-theoretic version) ensures that there is a unique homomorphism
$$
f_*:\Aut^0_{X/k}\lra \Aut^0_{Y/k}
$$
making $f$   equivariant with respect to $G=\Aut^0_{X/k}$.  The induced map of Lie algebras is likewise written  as 
$f_*:H^0(X,\Theta_{X/k})\ra H^0(Y,\Theta_{Y/k})$.
These maps are monomorphisms provided there is a $G$-stable schematically dense  open set
$U\subset X$ such that $f|U$ is an isomorphism to a schematically dense open set $V\subset Y$.
We then regard the maps as inclusions $\Aut^0_{X/k}\subset \Aut^0_{Y/k}$ and $H^0(X,\Theta_{X/k})\subset H^0(Y,\Theta_{Y/k})$.

Recall that a smooth proper  scheme $X$  is called \emph{variety of general type} provided $h^0(\O_X)=1$ and the dualizing sheaf
$\omega_X$ is a  \emph{big invertible sheaf}. Roughly speaking, this means that the function $t\mapsto h^0(\omega_X^{\otimes t})$ 
has a growths rate like a monomial of degree $n=\dim(X)$. According to  a   result of Martin-Deschamps and Lewin-M\'en\'egaux
\cite{Martin-Deschamps; Lewin-Menegaux 1978}, the group scheme $\Aut_{X/k}$ is then finite.

The  \emph{smooth surfaces of general type} $S$ can also be characterized by the condition that $(\omega_S\cdot D)>0$
for every movable curve $D\subset S$. Each such surface comes with a   \emph{canonical model} $Y$, where $\omega_Y$
is an ample invertible sheaf  and the singularities $\O_{Y,y}$ are \emph{rational double points}. Let us call these $Y$ \emph{canonical surfaces of general type}.
It is convenient to have a designation for the   modifications $X\ra Y$ whose  
singularities $\O_{X,x}$ are  rational double points; we propose to call them  \emph{RDP surfaces of general type}.
Special cases are the \emph{minimal surface of general type}, which arise from $S$ by successively contracting   $(-1)$-curves.
For more information, we refer to the monographs \cite{Barth et al 2004} and \cite{Badescu 2001}.
Let us record the following:

\begin{proposition}
\mylabel{prop:rdp surfaces general type automorphisms}
Let $X_k$ be a geometrically normal surface such that the minimal resolution of $X_{\bar{k}}=X\otimes_k\bar{k}$, where $\bar{k}$ the algebraic closure of $k$, is a smooth surface of general type. Then the group scheme $\Aut_{X/k}$ is finite.
\end{proposition}
Note that in particular the proposition applies to RDP surfaces of general type and hence they have finite automorphism group scheme.

\proof
Since $\Aut_{X/\bar{k}}\cong \mathrm{Aut}_{X/k}\otimes_k \bar{k}$, it suffices to prove the finiteness of $G=\mathrm{Aut}_{X/k}$ in the case when $k$ is algebraically closed.  Let $S\ra X$ be the minimal resolution of singularities. By the assumption, $S$ is a smooth surface of general type. Since $k$ is algebraically closed, $G_\red$ is a subgroup scheme of $G$. Since it is reduced, it fixes the singular locus of $X$ and therefore  there exists an injection $G_{red}\rightarrow \mathrm{Aut}_{U/k}$, where $U=X=\mathrm{Sing}(X)$, is the smooth locus of $X$. Hence $G_{red}$ is contained in the group of birational automorphisms $\mathrm{Bir}(S)$ of $S$, which is finite by~\cite{Martin-Deschamps; Lewin-Menegaux 1978}. Therefore $G_{\red}$ is also finite and hence $G$ is finite as well.

\qed

\medskip
From now on we assume $p>0$.
For each group scheme $G$  the Lie algebra $\lieg=\Lie(G)$ has, besides the bracket $[x,y]$ defined as commutator in the
ring of differential operators, also a \emph{$p$-map} $x^{[p]}$ stemming from $p$-fold composition of differential operators,
 and thus becomes a  \emph{restricted Lie algebra} (see \cite{Demazure; Gabriel 1970}, Chapter II, \S7 or 
\cite{Schroeer; Tziolas 2023}, Section 1 for the axioms). A vector $x\in\lieg$ is called \emph{$p$-closed} if the line $\idealh=kx$ is stable under the $p$-map,
in other words $x^{[p]}=\lambda x$ for some scalar $\lambda$, and thus becomes a commutative restricted subalgebra. 
In case $\lambda\neq 0$ the vector is called   \emph{multiplicative},
otherwise   \emph{additive}.
By the \emph{Demazure--Gabriel Correspondence} (\cite{Demazure; Gabriel 1970}, Chapter II, \S 7, Theorem 3.5)   
the functor $G\mapsto \lieg$ induces  an equivalence
between the categories of group schemes of finite type annihilated by the relative Frobenius  and the finite-dimensional
restricted Lie algebras. For convenience  such $G$ are called \emph{group schemes of height one}.
The group schemes $ \alpha_p^n=\GG_a^n[F]$ correspond to the standard vector space $\lieg=k^n$,
where both   bracket and $p$-map are trivial, while $\mu_p=\GG_m[p]$ corresponds to $\lieg=\liegl_1(k)$,
which is the one-dimensional standard vector space with $p$-map $\lambda^{[p]}=\lambda^p$.

Recall that   each associative  algebra $A$
becomes a  restricted Lie algebra, via $[x,y]=xy-yx$ and   $x^{[p]}=x^p$.
Write $\liegl_n(k)$ for the restricted Lie algebra arising from $A=\Mat_n(k)$.
It sits in two short exact sequences with kinks
\begin{equation}
\label{right four-term sequence}
\begin{tikzcd}[row sep=tiny, column sep=small ]
0\arrow[r]	& \liegl_1(k)\arrow[dr]\arrow[rr,dashed,"n"]	&  		& \liegl_1(k)\arrow[r] 		& 0\\		
		&						& \liegl_n(k)\arrow[ur]\arrow[dr]\\
0\arrow[r]	& \liesl_n(k)\arrow[ur]	\arrow[rr,dashed ]	&  		& \liepgl_n(k)\arrow[r]	& 0,\\
\end{tikzcd}
\end{equation}
where the inclusion of $\liegl_1(k)$ is given by the scalar matrices, and the surjection to $\liegl_1(k)$ is the trace map.
In turn, we get an exact sequence
$$
0\lra \Kernel(n|\liegl_1(k))\lra \liesl_n(k)\lra \liepgl_n(k)\lra \Cokernel(n|\liegl_1(k))\lra 0,
$$
where the outer terms are either  copies of $\liegl_1(k)$ or vanish, depending on the $p$-divisibility of $n$.
 
Throughout, we are particularly interested in the  restricted Lie algebra $\liesl_2(k)$, which by~\cite[Theorem 12.1]{Schroeer; Tziolas 2023} is one of the three possible restricted Lie algebras of the Frobenius kernel $\mathrm{Aut}_{X/k}[F]$ of a RDP surface $X$ of general type.

A basis of the restricted Lie algebra $\liesl_2(k)$ consists of the  matrices $h=(\begin{smallmatrix} 1&0\\0&-1\end{smallmatrix})$,  $e=(\begin{smallmatrix} 0&1\\0&0\end{smallmatrix})$ and 
$f=(\begin{smallmatrix} 0&0\\-1&0\end{smallmatrix})$. Bracket  and $p$-map  are determined by 
\begin{equation}
\label{eq:relations sl2}
[h,e]=2e,\quad [h,f]=-2f,\quad [e,f]=-h\quadand e^{[p]}=f^{[p]}=0,\quad h^{[p]}=h.
\end{equation} 
Note that the vector $h$ induces a subgroup  $\mu_p \hookrightarrow \mathrm{SL}_2$. This observation will be used later on.

For each traceless matrix $N=(\begin{smallmatrix}0&-d\\1&0\end{smallmatrix})$ in rational normal form one has 
$N^2=-dE$, which for $p\geq 3$ implies  $N^p=(-d)^{(p-1)/2}N$. Using that  every traceless matrix  
is similar to such an $N$, we obtain  
\begin{equation}
\label{eq:p-map in sl}
x^{[p]}=\begin{cases}
(\alpha^2-\beta\gamma)^{(p-1)/2}\cdot x	& \text{if $p\neq 2$;}\\
(\alpha^2-\beta\gamma)\cdot h				& \text{if $p=2$.}
\end{cases}
\end{equation} 
for each  vector $x=\alpha h+\beta e+\gamma f$ with coefficients $\alpha,\beta,\gamma\in k$. In turn, 
the homogeneous quadratic equation $\alpha^2-\beta\gamma=0$ defines the locus of additive vectors $x\in\liesl_2(k)$,
at least set-theoretically. 
The following two results are well-known, and we leave the verification to the reader.
 
\begin{proposition}
\mylabel{prop:sl2 for p odd}
For $p\neq 2$, the   canonical map $\liesl_2(k)\ra\liepgl_2(k)$ is bijective,     
every one-dimensional subspace $\lieh=kx$ is a restricted
subalgebra, and none of them is an ideal.
Moreover, there are no   two-dimensional ideals. 
\end{proposition}

This is in marked contrast to the situation at the remaining prime:

\begin{proposition}
\mylabel{prop:sl2 for p even}
For $p=2$ the kernel of $\liesl_2(k)\ra\liepgl_2(k)$ is generated  by $h$, the image  is a copy of $k^2$, and the resulting
extension $0\ra\liegl_1(k)\ra \liesl_2(k)\ra k^2\ra 0$ does not split. 
A   vector $x=\alpha h+\beta e+\gamma f$ generates a restricted subalgebra if and only if $\alpha^2-\beta\gamma=0$.
Furthermore, a non-zero subspace  $\idealh\subset\liesl_2(k)$ is a subalgebra if and only if it contains $h$,
and all such are restricted ideals. 
\end{proposition}

Note  that for $p\neq 2$, the height-one group scheme $\SL_2[F]$ is simple, while for $p=2$ it sits in  a non-split  
extension $1\ra \mu_2\ra \SL_2[F]\ra \alpha_2^{\oplus 2}\ra 0$.

For any restricted Lie algebra $\lieh$,  the non-zero homomorphisms $\liesl_2(k)\ra \lieh$   correspond to non-zero
triples $(e,h,f)$ from   $\lieh$ satisfying the   relations \eqref{eq:relations sl2}. 
Following Dynkin \cite{Dynkin 1957} and Bourbaki \cite{LIE 7-8} we call them
\emph{$\liesl_2$-triples}. The triple is called \emph{faithful} if the vectors $h,e,f\in \lieh$ are linearly independent.
Of course, the latter property is automatic in odd characteristics, by the preceding paragraph.
Combing this with the Demazure--Gabriel Correspondence, we see that for a given  scheme $X$
the $\liesl_2$-triples in $H^0(X,\Theta_{X/k})$ correspond to the non-trivial actions of the height-one group scheme   $\SL_2[F]$,
and the faithful triples become the faithful actions.

Let us say that a scheme $X$ \emph{admits an $\liesl_2$-triple} if this holds for the restricted Lie algebra $\lieh=H^0(X,\Theta_{X/k})$.
In \cite{Schroeer; Tziolas 2023} we raised the question whether or not there are surfaces of general type admitting faithful
$\liesl_2$-triples, and the goal of the present paper is to settle this problem.

Let us close this section with the following  observation:
Suppose $X$ is a     noetherian scheme with an $\liesl_2$-triple, with corresponding non-trivial action of 
$G=\SL_2[F]$. Recall that each non-zero vector $u=\alpha h+\beta e+\gamma f$  yields 
a subgroup scheme $H$ of order $p$. The resulting scheme of fixed points $X^H$ has the following regularity behaviour:

\begin{lemma}
\mylabel{lem:scheme of fixed points regular}
If the coordinates satisfy $\alpha^2-\beta\gamma\neq 0$, the scheme $X^H$ is regular at each point $x\in\Reg(X)$.
\end{lemma}

\proof
Without loss of generality we may assume that $X$ is the spectrum of a local noetherian ring that is regular.
According to \eqref{eq:p-map in sl}, the condition on the coordinates ensure that the group scheme $H$ is of multiplicative type,
and the assertion follows from \cite{Abramovich; Temkin 2018}, Proposition 5.1.16.
\qed

\medskip
This is  particularly useful if $X$ is regular  and  $X^H$ has a divisorial part, which then forms a   family of regular divisors
depending in an algebraic way on the line $ku\subset \liesl_2(k)$.
In one form or another, all what follows relies on this observation.

\section{Pencils and \texorpdfstring{$\liesl_2$}{sl2}-triples}
\mylabel{sec:pencils and triples}

Let $k$ be a ground field of characteristic $p>0$,
and $X$ be a proper scheme, together with a surjective morphism $f:X\ra\PP^1$ endowed with
a faithful $\liesl_2$-triple on $X$, and a compatible triple on $\PP^1$. In other words, we have 
a commutative diagram
\begin{equation}
\label{eq:triple in square}
\begin{CD}
\liesl_2(k)			@>>>	H^0(X,\Theta_{X/k})\\
@VVV					@VVV\\
H^0(\PP^1,\Theta_{\PP^1/k})	@>>>	H^0(X,f^*(\Theta_{\PP^1/k}))
\end{CD}
\end{equation} 
of restricted Lie algebras, where the upper horizontal map is injective, and the left vertical map is non-zero.
In geometric terms, the height-one group scheme $G=\SL_2[F]$ acts faithfully on $X$, and also in a compatible and non trivial way on the projective line.

Note that if $\O_{\PP^1}=f_*(\O_X)$, any $G$-action on $X$ automatically induces an action on $\PP^1$ by  Blanchard's Lemma.

Let  $\liesl_2(k)\cdot\O_X$ be the image of  the canonical map $\liesl_2(k)\otimes_k\O_X\ra \Theta_{X/k}$. Then $\liesl_2(k)\cdot \O_X$ is  a coherent subsheaf of the tangent sheaf that is stable 
under   bracket and $p$-map. 

Recall that an open set $U\subset \O_X$ containing the finite set $\Ass(\O_X)$ of all associated points
is termed \emph{schematically dense}.
Throughout this section, we assume that our datum satisfies the following two conditions:

\begin{enumerate}
\item The  $G$-action  on the projective line $\PP^1$  is fixed point free.
\item The coherent sheaf   $\liesl_2(k)\cdot\O_X$ is invertible on some schematically dense open set $U\subset X$.
\end{enumerate}
 
Note that if $p\geq 3$, condition $(i)$ is automatically satisfied. Indeed,  for $p\geq 3$,  $\mathfrak{sl}_2(k)$ is a simple restricted Lie algebra and since by assumption  the left vertical map in diagram $(\ref{eq:triple in square})$ is nontrivial, it must be injective,  and in fact  an isomorphism since both $\mathfrak{sl}_2(k)$ and $H^0(\PP^1,\Theta_{\PP^1/k})$ are three dimensional. Therefore, $G$ is isomorphic to $\mathrm{Aut}_{\mathbb{P^1}/k}[F]$ which is well known to act without fixed points on $\mathbb{P}^1_k$. 

The second condition holds if $X$ is a RDP surfaces of general type, and more generally for geometrically normal surfaces 
satisfying  $H^0(X,\omega_X^{\otimes-1})=0$, 
according to \cite{Schroeer; Tziolas 2023}, Corollary 6.6.
The combination of the above innocent assumptions has some remarkable consequences. Let us start with the following observation:

\begin{lemma}
\mylabel{lem:composite map bijective}
The  composite map $\liesl_2(k)\cdot\O_X\hookrightarrow \Theta_{X/k}\ra f^*(\Theta_{\PP^1/k})$ is bijective.
\end{lemma}

\proof
We will first show that the invertible sheaf $\Theta_{\PP^1/k}=\O_{\PP^1}(2)$ is globally generated by the image of the left vertical map in 
\eqref{eq:triple in square}. Taking into consideration that $\Theta_{\mathbb{P}^1/k}\cong \mathcal{O}_{\mathbb{P}^1}(2)$ is globally generated by two global sections without a common zero, it suffices to show that the image of the map in question contains two sections without a common zero.

Suppose that $p\geq 3$. In this case, as explained earlier after the statement of conditions $(i)$ and $(ii)$, the left vertical map is an isomorphism and hence its image contains two sections without a common zero.  

Suppose that $p=2$. In this case,  as explained in Proposition~\ref{prop:sl2 for p even},  $\mathfrak{sl}_2(k)$ is not simple and hence the left vertical arrow in $(\ref{eq:triple in square})$ may have a nontrivial kernel. Taking into consideration that the map in question is non zero, the dimension  of its image in $H^0(\PP^1,\Theta_{\PP^1/k})$ can be one, two or three. It cannot be one because if that was the case, the $G$-action on $\mathbb{P}^1$ would have a fixed point, which contradicts condition $(i)$. Hence the image is either two or three dimensional. Moreover, since we assume that the $G$-action on $\mathbb{P}^1$ is fixed point free, the sections of $H^0(\PP^1,\Theta_{\PP^1/k})$ generating the image must not have a common zero. 

Therefore, the invertible sheaf $\Theta_{\PP^1/k}=\O_{\PP^1}(2)$ is globally generated by the image of the vertical map on the left in 
\eqref{eq:triple in square}. This immediately implies that the map $\liesl_2(k)\cdot\O_X \ra f^*(\Theta_{\PP^1/k})$ is surjective and hence is an isomorphism on the schematically dense open set $U$ where by condition $(ii)$, $\liesl_2(k)\cdot\O_X$ is invertible. 

Suppose that the map $\liesl_2(k)\cdot\O_X \ra f^*(\Theta_{\PP^1/k})$ is not injective and that it has a nonzero kernel $\shN$. Let then $\zeta\in\Ass(\shN)$ be a non zero element. Since  the composite map 
\[
\liesl_2(k)\cdot\O_X\hookrightarrow \Theta_{X/k}\ra f^*(\Theta_{\PP^1/k})
\]
is an isomorphism on some schematically dense open set, we must have $\zeta\not\in\Ass(\O_X)$. Choose now a surjection $\O_{X,\zeta}^{\oplus r}\ra\Omega^1_{X/k,\zeta}$. The dual inclusion
$\Theta_{X/k,\zeta}\subset \O_{X,\zeta}^{\oplus r}$ shows  $\zeta\in \Ass(\O_X)$,  contradiction. Hence $\shN=0$.
\qed

\medskip
In other words, the canonical map $ \Theta_{X/k}\ra f^*(\Theta_{\PP^1/k})$ of coherent sheaves is split surjective,
and the foliation $\liesl_2(k)\cdot\O_X$ defines a splitting. Dualizing the  exact sequence
$f^*(\Omega^1_{\PP^1/k})\ra\Omega^1_{X/k}\ra \Omega^1_{X/\PP^1}\ra 0$ thus induces
a split extension coming with a decomposition
\begin{equation}
\label{eq:splitting tangent sheaf}
\Theta_{X/k}=\Theta_{X/\PP^1}\oplus f^*(\Theta_{\PP^1/k}),
\end{equation} 
where the second summand coincides with $\liesl_2(k)\cdot\O_X$
Informally speaking,  the   group scheme $\Aut_{X/k}$ looks  infinitesimally like a product. There are two useful consequences of this. 

The first one is that  $\liesl_2(k)\cdot \O_X$ is a foliation in the usual sense. This means that the quotient $\Theta_{X/k}/\mathfrak{sl}_2(k)\cdot \mathcal{O}_X$ is torsion free, a property that may not hold in general without assumptions $(i)$ and $(ii)$. This immediately follows from the statement of Lemma~\ref{lem:composite map bijective} and the equation $(\ref{eq:splitting tangent sheaf})$.

The second consequence is the following.

\begin{proposition}
\mylabel{prop:action faithful and p odd}
The induced  $G$-action   on $\PP^1$ is faithful, and the characteristic must be $p\geq 3$.
\end{proposition}

\proof
In the commutative diagram \eqref{eq:triple in square}, the upper vertical map is injective by our standing assumption.
In light of the splitting \eqref{eq:splitting tangent sheaf}, the composite mapping  $\liesl_2(k)\ra H^0(X,f^*(\Theta_{\PP^1/k}))$ remains injective,
and it follows that vertical map on the left is injective as well. This means that the $G$-action on 
the projective line is faithful. 
For  $p=2$, the three-dimensional Lie algebras $\liesl_2(k)$ and $\liepgl_2(k)$ are not isomorphic, which 
produces a contradiction.
\qed

\medskip
Let $B=\Spec(\O_X)$ be the Stein factorization for $f:X\ra \PP^1$. By Blanchard's Lemma this carries
another compatible $\liesl_2$-triple.

\begin{proposition}
\mylabel{prop:stein factorization projective line}
If $X$ is geometrically normal and geometrically connected, then $B$ is a twisted form of the projective line.
\end{proposition}

\proof
It suffices to treat the case that $k$ is algebraically closed. As $X$ is normal and connected, the same holds for $B$.
Since curves of genus $g\geq 2$, and elliptic curves alike, do not admit $\liesl_2$-triples, the only possibility is 
$B\simeq\PP^1$.
\qed

\medskip
The following will  play a key role throughout:

\begin{proposition}
\mylabel{prop:morphism is smooth}
If the proper scheme $X$ is smooth then  the     morphism $f:X\ra\PP^1$ is smooth as well.
\end{proposition}

\proof
Without loss of generality we may assume that the ground field $k$ is algebraically closed.
For each generic point $\eta\in X$, the corresponding irreducible component is a connected component $U\subset X$,
and its image $f(U)\subset \PP^1$ is $G$-stable. Since the action on the projective line is fixed point free,
$f(\eta)\in\PP^1$ must be the generic point, hence the morphism $f$ is flat.

The sheaf of K\"ahler differentials $\Omega^1_{X/k}$ is locally free of rank $n=\dim(X)$. 
In the exact sequence
$$
f^*(\Omega^1_{\PP^1/k})\stackrel{\varphi}{\lra}\Omega^1_{X/k}\lra \Omega^1_{X/\PP^1}\lra 0,
$$
the terms on the left are locally free. The dual map $\varphi^\vee$ is split surjective,
therefore   $\varphi=\varphi^{\vee\vee}$ is split injective. It follows that $\Omega^1_{X/\PP^1}$ 
is locally free of rank $n-1$, hence the fibers of $f$ are smooth.
Summing up, the morphism $f:X\ra \PP^1$ if flat with smooth fibers, hence smooth.
\qed

\medskip
Let us record the following consequence:

\begin{corollary}
\mylabel{cor:family of varieties of general type}
Suppose that the scheme  $X$ is smooth of dimension $m+1$, that the morphism $f:X\ra\PP^1$ is in Stein factorization,
and that the dualizing sheaf $\omega_X$ restricts to a  big invertible sheaf 
on the generic fiber $X_\eta$. Then $X$ is a family of $m$-dimensional smooth varieties of general type
parametrized by $\PP^1$.
\end{corollary}

\proof
It suffices to treat the case that the ground field $k$ is algebraically closed.
For each closed point $t\in \PP^1$, the Adjunction Formula gives $\omega_{X_t}=\omega_X|X_t$.
So by the  Semicontinuity Theorem, the dualizing sheaf $\omega_{X_t}$ is big.
The fiber $X_t$ is smooth, according to the Theorem, and geometrically connected by  $\O_{\PP^1}=f_*(\O_X)$.
Thus $X$ is a family of smooth varieties of general type, obviously of dimension $m=\dim(X)-1$.
\qed

\section{The inertia map}
\mylabel{sec:inertia map}

Let $k$ be a ground field of characteristic $p>0$, and $Y$ be a reduced scheme of finite type endowed with
an $\liesl_2$-triple, in other words, a non-trivial action of the height-one group scheme $G=\SL_2[F]$.
We   also write  $\lieg=\liesl_2(k)$ for the  Lie algebra.
The  cartesian diagram
$$
\begin{CD}
I	@>>>	Y\\
@VVV		@VV\Delta V\\
G\times Y	@>>>	Y\times Y,
\end{CD}
$$
where the lower map is given by $(g,y)\mapsto (gy,y)$, defines the \emph{inertia} $I$,
which is a relative  group scheme over $Y$ whose structure morphism is finite, but usually far from flat.
The  fibers $I_y\subset G\otimes\kappa(y)$ can be seen are the  \emph{stabilizers} $G_y$
with respect to $y$ viewed as rational point on the base-change $Y\otimes\kappa(y)$, and correspond to the restricted  subalgebra
$$
\Lie(G_y)\subset \lieg\otimes\kappa(y) =\liesl_2(\kappa(y)).
$$
Let $Y_d\subset Y$, $0\leq d\leq 3$
be the  subscheme of points $y\in Y$ where the order the stabilizer takes the value $|I_y|=p^d$, The stratum  $Y_0$ and can be seen as the locus where the $G$-action is free,
whereas $Y_3$  equals the \emph{scheme of fixed points} $Y^G$. This yields a stratification in the sense that $Y_{i+1}\subset \overline{Y}_i\smallsetminus Y_i$. By upper semicontinuity of the length of fibers of the finite morphism $I\rightarrow Y$, $Y_0$ is open in $Y$, $Y_3$ is closed  in $Y$ and $Y_{i+1}$ is open in $ \overline{Y}_i\smallsetminus Y_i$, for $i=0,1$.

\emph{In what follows, we assume that   $Y_0=Y_1=\varnothing$,    and furthermore that $Y_2$ is dense.} Then 
our stratification simplifies to  the closed
scheme $Y^G=Y_3$ and the complementary open set $Y\smallsetminus Y^G=Y_2$. We see that $Y^G=\varnothing$
if and only if all the Lie algebras $\Lie(G_y)$, $y\in Y$ are two-dimensional. In this situation, the $\liesl_2$-triple is called
\emph{fixed point free}. In any case, the inertia $I\ra Y$ becomes flat   outside the closed set  $Y^G$,
and we obtain a  family of height one group schemes of length $p^2$. 

By the Demazure-Gabriel correspondence, such families can be parametrized by parametrizing their corresponding restricted Lie algebras. This can be achieved as follows. Over $U=Y\smallsetminus Y^G$, there exists a short exact sequence
\[
0 \rightarrow N \rightarrow \mathfrak{g}\otimes \mathcal{O}_U \stackrel{\theta}{\rightarrow} \mathfrak{g}\cdot \mathcal{O}_U \rightarrow 0.
\]
At any point $y\in U$, the kernel of the map $\theta \otimes k(y)$ is the restricted Lie algebra $\mathrm{Lie}(G_y)\subset \mathfrak{g} \otimes k(y)$, where $G_y$ is the stabilizer group scheme of $y\in U$. Considering that $\mathrm{Lie}(G_y)$ is 2-dimensional and $\mathfrak{g}$ is 3-dimensional, it immediately follows that 
\[
\dim_{k(y)} ( \mathfrak{g}\cdot \mathcal{O}_U)\otimes k(y)=1,
\]
for any $y\in U$. By taking into consideration that since $ \mathfrak{g}\cdot \mathcal{O}_U \subset \Theta_{U/k}$,  $\mathfrak{g}\cdot \mathcal{O}_U$ is torsion free, this implies that $ \mathfrak{g}\cdot \mathcal{O}_U$ is invertible and hence $N$ is locally free of rank $2$.  One may now utilise the theory of the Grassmanian varieties and conclude that there exists a morphism 
$$
f_\inert:Y\smallsetminus Y^G\lra  \Grass^2(\lieg),\quad y\longmapsto [\Lie(G_y)].
$$
into the Grassmann variety of two-dimensional vector subspaces of $\mathfrak{g}$,  where $[\Lie(G_y)]\in \Grass^2(\lieg)$ is the closed point of  the Grassmann variety corresponding to $\Lie(G_y)$. Moreover, from the universal property of the Grassmann variety, $\mathfrak{g}\cdot \mathcal{O}_U=f_\inert^{\ast}L$, where $L$ is the universal quotient invertible sheaf of $\Grass^2(\lieg)$. The map  $f_\inert$ is of fundamental importance, and will be called the \emph{inertia map}. 

The inertia map is not constant, since our action is faithful.  
Neither can it be surjective,  because there are subspaces that fail to be subalgebras.
To determine the image  we must find the conditions under which a 2-dimensional subspace of $\mathfrak{g}$ is a subalgebra of $\mathfrak{g}$.  In order to do this we use the Pl\"ucker coordinates of a 2-dimensional subspace of $\mathfrak{g}$. 

Recall from the classical theory of the Grassmannian variety, that there exists an embedding $\Grass^2(\lieg)\subset \mathbb{P}(\Lambda^2\lieg)$, which in this particular case is an equality for dimension reasons. This embedding is called the  Pl\"ucker  embedding and  is defined as follows. The $k$-basis $h, e, f$ of $\mathfrak{g}$ defines the basis $h \wedge f$, $h\wedge e$ and $e \wedge f$ of $\Lambda^2\lieg$. Let $L$ be a two dimensional subspace of $\mathfrak{g}$ defined by the two linearly independent vectors $v=\alpha h+\beta e+\gamma f$ and $v'=\alpha' h+\beta'e+\gamma' f$. The subspace $L$ corresponds to the point 
\[
v\wedge v^{\prime}=(\alpha\gamma'-\gamma\alpha' )h \wedge f+(\alpha\beta'-\beta\alpha' )h\wedge e +(\beta\gamma'-\gamma\beta')e \wedge f 
\] 
of $\mathbb{P}(\Lambda^2\lieg)$. The numbers $T_0=\alpha\gamma'-\gamma\alpha' $, $T_1=\alpha\beta'-\beta\alpha'$ and $T_2= \beta\gamma'-\gamma\beta'$
are called the Pl\"ucker coordinates of the subspace. 

Given linearly independent vectors $v=\alpha h+\beta e+\gamma f$ and $v'=\alpha' h+\beta'e+\gamma' f$,
the 2-dimensional subspace  $kv+kv'$ of $\mathfrak{g}$ is a subalgebra of $\mathfrak{g}$ if and only if $v\wedge v'\wedge [v,v']=0$ in $\Lambda^3\mathfrak{g}$. A  straightforward computation in the Grassmann algebra  using the relations \eqref{eq:relations sl2} in our Lie algebra $\lieg=\liesl_2(k)$ shows
$$
v\wedge v'\wedge [v,v']=  (T_0^2-4T_1T_2)\cdot ( h\wedge e \wedge f),  
$$
where  $T_0$, $T_1$ and $T_2$
are the Pl\"ucker coordinates of the subspace. Considering that in  our situation, the Pl\"ucker embedding $\Grass^2(\lieg)\subset \mathbb{P}(\Lambda^2\lieg)$ is an equality,  the above computation reveals:

\begin{proposition}
\mylabel{prop:image inertia map}
Inside  $\mathbb{P}(\Lambda^2\lieg)$, the scheme of lines stemming from  subalgebras $\idealh\subset \lieg$ is the curve
defined by the quadratic equation $T_0^2-4T_1T_2=0$.
\end{proposition}
Let us  write $B=V_+(T_0^2-4T_1T_2)$   and call it the \emph{curve of subalgebras}. If $p\not=2$ then this is a smooth conic in $\mathbb{P}^2$ and hence isomorphic to $\mathbb{P}^1$. If on the other hand $p=2$, this is a double line in $\mathbb{P}^2$, in particular it is non reduced. In any case, the inertia map factors as 
\begin{equation}
\label{eq:inertia map factorization}
f_\inert:Y\smallsetminus Y^G\lra  B=V_+(T_0^2-4T_1T_2)\subset \mathbb{P}(\Lambda^2\lieg).
\end{equation} 

On the Grassmann variety $V=\Grass^2(\lieg)= \mathbb{P}(\Lambda^2\lieg)$, the universal quotient invertible sheaf of  $\Grass^2(\lieg)$ is $\mathcal{O}_{V}(1)$ and hence 
$\mathfrak{g}\cdot \mathcal{O}_U=f_\inert^{\ast}\mathcal{O}_V(1)$. This implies that $\mathfrak{g}\cdot \mathcal{O}_U$ is globally generated by the pullbacks of the three global sections of $\mathcal{O}_V(1)$ which generate it globally and  which correspond to the images of $e,f, h$ under the natural isomorphism $\mathfrak{g} \rightarrow H^0(V,\mathcal{O}_V(1))$. Therefore there exists a canonical injection
\begin{equation}
\label{eq:new interpretation of triple}
\liesl_2(k)=\lieg = H^0(V,\O_V(1)) \stackrel{f_\inert^*}{\lra}  H^0(U,\O_U(1))\subset H^0(U,\Theta_{X/k})
\end{equation} 
which yields the restriction of the original $\liesl_2$-triple on $Y$. This has a surprising consequence:

\begin{proposition}
\mylabel{prop:char not two}
In the above situation, the characteristic must be  $p\geq 3$.
\end{proposition}

\proof
Suppose  $p=2$.
Then the  curve  of   subalgebras   becomes  $V_+(T_1^2)$.
Since our scheme $Y$ is reduced, the inertia map factors over the line $L=V_+(T_1)$.
It follows that the injection \eqref{eq:new interpretation of triple} factors over the two-dimensional vector space $H^0(L,\O_L(1))$, contradiction.
\qed

\medskip
\emph{So from this point onward we know $p\neq 2$.}
Our group scheme $G=\SL_2[F]$ acts on the scheme $Y$,   on itself by conjugacy, and on the Lie algebra via the adjoint map 
$\Ad: G\ra\GL(\lieg)$. 
Since $G$ has height-one, the latter is already determined by its derivative
$$
\ad:\lieg\lra\liegl(\lieg),\quad u\longmapsto (v\mapsto [u,v]).
$$
To make this explicit, consider $v=\alpha h+\beta e+\gamma f$ and $v'=\alpha' h+\beta'e+\gamma' f$ as above.
The adjoint action of $u=\pi h + \varphi e+\psi f$ on the plane $kv+kv'$ takes the form
$v\wedge v'\mapsto [u,v]\wedge[u,v']$. A straightforward computation shows that this has the matrix interpretation
\begin{equation}
\label{eq:induced action matrix interpration}
u\longmapsto \begin{pmatrix}
-4\pi^2	& -2 \pi\varphi	&  2\pi\varphi\\
 4\pi\varphi&  2\varphi\psi	& -2\varphi^2\\
-4\pi\psi	& -2\psi^2		&  2\varphi\psi
\end{pmatrix}
\end{equation}
with respect to the basis $h\wedge e$, $h\wedge f$ and $h \wedge e$ of $\Lambda^2 \mathfrak{g}$. This corresponds to an automorphism of $\mathbb{P}(\Lambda^2\mathfrak{g})$ defining the $G$-action on $\mathbb{P}(\Lambda^2\mathfrak{g})$.
One easily checks directly that this $G$-action stabilizes the closed subscheme defined by  $T_0^2-4T_1T_2=0$.
%

\begin{proposition}
\mylabel{prop:induced action fixed point free}
The  inertia map $f_\inert:Y\smallsetminus Y^G\ra \mathbb{P}(\Lambda^2\lieg)$ is $G$-equivariant, and the  induced action on $B=V_+(T_0^2-4T_1T_2)$
is fixed point free.
\end{proposition}

\proof
Equivariance of the inertia map can be checked for the functor of points, 
and then boils down to the set-theoretical fact that if a group $G$ acts on a set $Y$, then  $\sigma G_y\sigma^{-1}=G_{\sigma y}$, for any $\sigma \in G$ and $y\in Y$.
From \eqref{eq:induced action matrix interpration} we see that $u=e$ 
has the effect $(1:2:1)\mapsto (2:0:0)$. Both of which belong to the curve of subalgebras, and it
follows that the $G$-action is non-trivial. Since the group scheme $G$ is simple, the action must be faithful.
This identifies $G$ with the Frobenius kernel of the automorphism group scheme for $V_+(T_0^2-4T_1T_2)=\PP^1$,
and the latter acts without fixed points.
\qed

\medskip
Properness of the curve of subalgebras  has the following consequence:

\begin{proposition}
\mylabel{prop:codimension one fixed point free}
No point  $\zeta\in Y$ where the local ring $R=\O_{Y,\zeta}$ is regular and one-dimensional belongs to the scheme of fixed points $Y^G$.
\end{proposition}

\proof
By our standing assumption, the open set $Y\smallsetminus Y^G$ is dense, so the inertia map is defined at the generic point $\eta\in Y$
corresponding to $F=\Frac(R)$.
By the Valuative Criterion, the morphism $f_\inert:\Spec(F)\ra \Grass^2(\lieg)$ extends
to $\Spec(R)$. Hence the inertia map is defined on some open set $U$ which contains both $\zeta$ and $U_0=Y\smallsetminus Y^G$.
The resulting morphism $U\ra\PP^1$ is equivariant, because this holds on a schematically dense open set $U_0$.
All    open sets are stable with respect to the infinitesimal group scheme $G$.
Since the $G$-action on $V_+(T_0^2-4T_1T_2)$ is fixed point free, the same holds on $U$, and therefore $\zeta\not\in Y^G$.
\qed

\medskip
So if  $Y$ is normal, the closed subscheme $Y^G$ has codimension at least two. 
We now regard \eqref{eq:inertia map factorization} as   rational map $f_\inert:Y\dashrightarrow B$
and seek to resolve the locus of indeterminacy. To this end consider the graph
$$
\Gamma_\inert\subset (Y\smallsetminus Y^G)\times B,
$$
and write $\tilde{Y}$ for its schematic closure inside $Y\times B$.
This closed subscheme is $G$-stable, according to \cite{Brion; Schroeer 2023}, Lemma 2.1. 
By construction, it is   an \emph{equivariant compactification}
of $Y\smallsetminus Y^G$,    and the projection $r:\tilde{Y}\ra Y$ is an \emph{equivariant modification}.
Note that $\tilde{Y}$ might be non-normal, even if $Y$ is normal.
The resulting inertia map 
$$
\tilde{f}_\inert:\tilde{Y}\lra B=  V_+(T_0^2-4T_1T_2)\subset\Grass^2(\Lambda^2\lieg)
$$
coincides with the composition $f_\inert\circ r$.  Our second main result is that this construction
``resolves  the fixed points'' for the action of $G=\SL_2[F]$:

\begin{theorem}
\mylabel{thm:resolution of fixed points}
The $G$-action on $\tilde{Y}$ is fixed point free.  Moreover, the coherent subsheaf  $\lieg\cdot\O_\tY\subset\Theta_{\tY/k}$
is invertible, the  canonical map $\Theta_{\tY/k}\ra \tf^*_\inert(\Theta_{B/k})$
is surjective, and  the foliation   $\lieg\cdot\O_\tY$ yields a splitting.
\end{theorem}

\proof
The first statement simply follows from the fact that the $G$-action on $B$ is fixed point free.
The remaining assertions are consequences of Lemma \ref{lem:composite map bijective}.
\qed

\medskip
In light of the above, we call    $r:\tilde{Y}\ra Y$ the \emph{minimal resolution of fixed points}.
It is quite remarkable that this exists.


\section{Surfaces of general type}
\mylabel{sec:surfaces}
 
Let $k$ be a ground field of characteristic $p>0$. The following result answers a question 
raised in \cite{Schroeer; Tziolas 2023}, and was our initial motivation for this research:

\begin{theorem}
\mylabel{thm:general type no triple}
No proper smooth surface of general type has a faithful $\liesl_2$-triple.
\end{theorem}

The proof requires some preparation, and will be given at the end of this section. Our methods actually give insight for the geometry of  $\liesl_2$-triples on 
other surfaces as well, and we consider the following general setting: Let $Y$ be a   normal proper surface with $h^0(\O_Y)=1$,
endowed with an $\liesl_2$-triple such that the subsheaf $\liesl_2(k)\cdot\O_Y\subset\Theta_{Y/k}$ has rank one.
In other words, we have a non-trivial action of the  height-one group scheme $G=\SL_2[F]$, with Lie algebra $\lieg=\liesl_2(k)$,
such that for the generic stabilizer $G_\eta$ has two-dimensional  Lie algebra.
Consequently, the theory developed in Section \ref{sec:inertia map} applies in our situation:
$\Lie(G_y)$ is two-dimensional for all $y\not\in Y^G$, the  
scheme of fixed points $Y^G$ is finite (Proposition \ref{prop:codimension one fixed point free}), and
the characteristic  is different from two (Proposition \ref{prop:char not two}).

Let  $\tY\ra Y$ be the minimal resolution of fixed points.  
Since the group scheme $G$ is infinitesimal, the  action on $Y$ does not necessarily extend
to the normalization, so $\tY$ can be non-normal. However, the situation improves  under suitable  assumption:

\begin{proposition}
\mylabel{prop:equivariant normalization}
Let $X\ra \tY$ be the normalization.
If the surface  $Y$  has only rational singularities, the $G$-action extends to $X$, and the latter
has only rational singularities.
\end{proposition}

\proof
We recursively define a   sequence
$$
Y=Y_0\longleftarrow Y_1\longleftarrow Y_3 \longleftarrow \ldots
$$
of proper normal $G$-surface $Y_i$ and $G$-equivariant  modifications as follows:
To start with set $Y_0=Y$. Suppose now that $Y_i$ is already defined.
If the $G$-action is fixed point free we set $Y_{i+1}=Y_i$. Otherwise   pick some $z_i\in Y_i^G$, regard $Z_i=\{z_i\}$
as reduced closed subscheme, 
and define $Y_{i+1}=\Bl_{Z_i}(Y_i)$. Then the $G$-action on $Y_i$ extends to $Y_{i+1}$, 
and the canonical morphism $Y_{i+1}\ra Y_i$ is equivariant. According to \cite{Lipman 1969}, Proposition 8.1 the surface $Y_{i+1}$
stays normal, and has only rational singularities. This concludes the recursive definition. 

According to \cite{Lipman 1969}, Theorem 26.1 the rational map $f_\inert:Y\dashrightarrow \PP^1$
becomes a morphism on     $Y_r$ for some index $r\geq 0$, so the $G$-action on $Y_r$ is fixed point free.
Therefore  the canonical morphism $Y_r\ra Y$
factors over $\tilde{Y}$,   thus also over the normalization $X$. 
The induced morphism $h:X\ra Y$ is in Stein factorization, so by Blanchard's Lemma there is a $G$-action on $Y$ for which
$h$ is equivariant. It follows that the normalization map $X\ra \tY$ is equivariant. 

To see that $X$ has only rational singularities, choose a resolution of singularities $X'\ra Y_r$.
The Leray--Serre spectral sequence for the composition of $g:X'\ra Y$ and $h:X\ra Y$ gives an exact sequence
$$
R^1(h\circ g)_*(\O_{X'})\lra h_* (R^1g_*(\O_{X'}))\lra R^2h_*(\O_Y).
$$
The term on the left vanishes because the singularities on $X$ are rational, and the term on the right is trivial
by the Theorem of Formal Functions. Consequently $R^1g_*(\O_{X'})=0$.
\qed

\medskip
Note that the $G$-action on $X$ is fixed-point free, since this already holds for $\tilde{Y}$. 
Applying Theorem \ref{thm:resolution of fixed points} with $X$ instead of $Y$, we see that $\lieg\cdot\O_X\subset\Theta_{X/k}$
is invertible, the canonical map $\Theta_{X/k}\ra \tf^*_\inert(\Theta_{B/k})$
is surjective, and the foliation $\lieg\cdot\O_X$ yields a splitting.

\begin{proposition}
\mylabel{prop:modifications with good properties}
Suppose $Y$ is regular, or has only rational double points, or has only rational singularities.
Then there is a $G$-modification $X\ra Y$ with some normal surface $X$ that has the respective properties,
and where the $G$-action is fixed point free.
\end{proposition}

\proof
Use the sequence $Y=Y_0\leftarrow Y_1\leftarrow\ldots\leftarrow Y_r$ of blowing-ups as in the preceding proof.
The $Y_i$ are regular, or have only rational double points, or   have only rational singularities if the corresponding
property holds for $Y$. We thus may set $X=Y_r$. 
\qed

\begin{proposition}
\mylabel{prop:modifications is  family of curves}
Suppose that the ground field $k$ is algebraically closed and $Y$ is regular. Then one of the following holds.
\begin{enumerate}
\item $Y$ is isomorphic to a product  $C \times \mathbb{P}^1$, where $C$ is a smooth proper curve of genus $g\geq 1$.
\item $Y$ is isomorphic to the Hirzebruch  surface $\mathbf{F}_{pn}=\mathbb{P}(\mathcal{E})$, where $\mathcal{E}\cong \mathcal{O}_{\mathbb{P}^1}\oplus \mathcal{O}_{\mathbb{P}^1}(-pn)$, for some $n\geq 0$. Moreover, $Y$ is the base-change  of the Hirzebruch surface $\mathbf{F}_n$, with respect to the relative Frobenius map $F:\PP^1\ra\PP^1$ on the base.
\end{enumerate}
In all cases  there exists a unique, up to automorphisms of $\mathfrak{sl}_2(k)$, $\mathfrak{sl}_2$-triple on $Y$, except in the case when $X\cong \mathbb{P}^1\times \mathbb{P}^1$, when there are two.
\end{proposition}

\begin{proof}

According to Proposition \ref{prop:modifications with good properties}, there is a $G$-modification $\phi \colon X\ra Y$ with some regular surface $X$ where the $G$-action is fixed point free.
By Proposition \ref{prop:induced action fixed point free}, the inertia map $f_\inert$ is equivariant, and the induced action on $B$ is fixed-point free. 
Since $k$ is algebraically closed,  the Stein factorization takes the form $f:X\ra\PP^1$, and this morphism is smooth 
(Proposition \ref{prop:stein factorization projective line} and
Proposition \ref{prop:morphism is smooth}). In turn, we may view $X$ as a family of   curves, of a certain genus $g\geq 0$.

Suppose now that $g\geq 2$.
Let $\PP^1\ra\shM_g$ be the classifying map to the Deligne--Mumford stack.
Such maps have zero-dimensional image: This was first established by Parshin \cite{Parshin 1968} 
and Arakelov \cite{Arakelov 1971} in characteristic zero, and later extended to positive characteristics
by Szpiro (\cite{Szpiro 1979}, Theorem 3.3).
It follows that our  family of curves
becomes constant on some connected  \'etale covering of the projective line. Since $k$ is separably closed,
the only connected such covering is the identity.  We thus have $X=C\times \PP^1$ for some curve $C$.
Nguyen's further generalization (\cite{Nguyen 1998}, Theorem 2) yields the product structure for  the case $g=1$. 

Note also that a surface isomorphic to a product $C\times \mathbb{P}^1$,  does not contain any $(-1)$-curves. Therefore the map $\phi\colon X\ra Y$ is an isomorphism and hence actually $Y= C \times \mathbb{P}^1$. This shows part $(i)$ of Proposition~\ref{prop:modifications is  family of curves}. Considering now that 
\[
H^0(Y,\Theta_{Y/k})=H^0(C,\Theta_{C/k})\oplus H^0(\mathbb{P}^1,\Theta_{\mathbb{P}^1/k})
\]
and that if $g\geq 1$,  $C$ does not have any $\mathfrak{sl}_2$-triples, since $H^0(C,\Theta_{C/k})$ is either zero or $k$, it follows that $Y$ has exactly one, up to automorphisms of $\mathfrak{sl}_2(k)$,  $\mathfrak{sl}_2$-triple coming from $\mathbb{P}^1$.

Suppose now that $g=0$. Then $f\colon X \rightarrow \mathbb{P}^1$ is a ruled surface. Hence either the Hirzebruch surface $X=\mathbf{F}_m$, $m\geq 1$, or $X=\mathbb{P}^1\times \mathbb{P}^1$. $X=\mathbb{P}^1\times \mathbb{P}^1$ has two $\mathfrak{sl}_2$-triples coming from the two factors. Suppose that $X=\mathbf{F}_m$, $m\geq 1$. By Lemma~\ref{lem:composite map bijective}, the exact sequence
\begin{equation}\label{eq 1000}
0 \rightarrow \Theta_{X/\mathbb{P}^1} \rightarrow \Theta_{X/k} \rightarrow f^{\ast}\Theta_{\mathbb{P}^1/k}\rightarrow 0
\end{equation}
is split exact. By~\cite[Theorem 4.2]{Ganong; Russell 1982}, this happens exactly when the characteristic $p$ divides $m$, i.e., $m=pn$, for some $n>0$. Let $\sigma \colon f^{\ast}\Theta_{\mathbb{P}^1/k}\rightarrow \Theta_{Y/k}$ be  a splitting 
of~(\ref{eq 1000}). Then $\sigma ( f^{\ast}\Theta_{\mathbb{P}^1/k})\subset \Theta_{X/k}$ defines an $\mathfrak{sl}_2$-triple on $X.$ By taking into consideration that every $\mathfrak{sl}_2$-triple defines a ruling $X\rightarrow \mathbb{P}^1$ by use of the inertia map and that a Hirzebruch surface $\mathbf{F}_m$ with $m\geq 1$, has exactly one ruling, there exists a unique, up to isomorphism, $\mathfrak{sl}_2$-triple on $X$.

Finally, since $\mathbf{F}_{pn}$ is minimal, it follows that $X$ does not contain any $(-1)$ curves. Hence $\phi \colon X\rightarrow Y$ is an isomorphism and therefore $Y\cong \mathbf{F}_{pn}$. This concludes the proof of 
Proposition~\ref{prop:modifications is  family of curves}.
\end{proof}

\emph{Proof for Theorem \ref{thm:general type no triple}.}
Without loss of generality we may assume that $k$ is algebraically closed.
Seeking a contradiction, we assume that there is a smooth surface $X$ of general type having a faithful $\liesl_2$-triple.
Since the $G$-action is non-trivial, the corresponding foliation $\lieg\cdot\O_X\subset\Theta_{X/k}$ is non-zero.
The rank of this coherent sheaf is at least one, and bounded above by the so-called \emph{foliation rank} $r\geq 0$,
an invariant introduced in \cite{Schroeer; Tziolas 2023}, Section 6.
The condition $h^0(\omega_X^{\otimes-1})=0$ ensures $r\leq 1$, according to loc.\ cit.,  Corollary 6.6.
So for the function field $F=k(X)$, the kernel of the canonical map $\lieg\otimes_kF\ra \lieg\cdot F\subset\Theta_{X/k,\eta}$
has dimension $d=2$. By \cite{Schroeer; Tziolas 2023}, Proposition 5.3 the Lie algebra of the generic stabilizer $G_\eta$
has dimension $d=2$. 

Thus the above results apply, and by Proposition \ref{prop:modifications is  family of curves} we may furthermore assume that 
$X$ is a family of smooth curves  over   the projective line $\PP^1$, which by 
Corollary \ref{cor:family of varieties of general type}    have genus $g\geq 2$, and
thus   $X=C\times\PP^1$.
Let $D=\{c\}\times\PP^1$ be some closed fiber with respect to the first projection. This curve is movable and $\omega_X$ is big, hence $(\omega_X\cdot D)>0$.
On the other hand, the Adjunction Formula gives $(\omega_X\cdot D)=-2$, contradiction.
\qed

\section{Lefschetz pencils and  alterations}
\mylabel{sec:lefschetz pencils}

Let $k$ be a ground field of characteristic $p>0$.
The goal of this section is to construct
RDP surfaces  of general type admitting   $\liesl_2$-triples. Note that by  Blanchard's Lemma, this also produces
$\liesl_2$-triples on their canonical model.
In what follows, we work with the height-one group scheme $\PGL_2[F]$ instead of $\SL_2[F]$,
in order to keep $p=2$ included.   Our main result in this direction is:

\begin{theorem}
\mylabel{thm:existence with higher height}
For each smooth connected proper surface $S$    and each integer $n\geq 0$, 
there is an  alteration $X\ra S$ with purely inseparable function field extension of $\deg(X/S)=p^n$ such that the following holds:
\begin{enumerate}
\item There is a  faithful action of the group scheme $\PGL_2[F^n]$ on $X$.
\item The surface $X$ is normal, the  singularities are rational double points  of type $A_l$ with $l=p^n-1$,  
and  the tangent sheaf  $\Theta_{X/k}$ is locally free.
\item If $S$ is a smooth surface of general type, one may choose $X$ as a RDP surface of general type.
\end{enumerate}
\end{theorem}

The proof is constructive in nature and will be given in the course of this section,
after some preparatory considerations.  

Let $\shL$ be a very ample invertible sheaf on our smooth surface $S$,
and  $E\subset H^0(S,\shL)$ be some 2-dimensional vector subspace. Choose a basis $s_0,s_1\in E$ and consider
the resulting \emph{pencil} of curves $C_t\subset S$ parameterized by $t\in\PP^1$, and denote the  \emph{axis} by
$Z=\bigcap C_t=C_0\cap C_\infty$.
It  is called   \emph{Lefschetz pencil} if almost all $C_t$ are smooth, each singular $C_t$ has but one singularity
that is furthermore \'etale locally given by the equation $xy=0$, and the axis $Z$ is finite and smooth.
Upon blowing-up with the  axis $Z$ as center, one obtains a fibration $\tilde{f}\colon \tilde{S}=\Bl_Z(S)\ra\PP^1$. 
The fibration $\tilde{f}\colon \tilde{S} \rightarrow \mathbb{P}^1$ is called  the \emph{Lefschetz fibration}.
According to \cite{SGA 7b}, Expos\'e XVII, Theorem 2.5 the Lefschetz pencils form a dense open set inside the Grassmann variety of planes,
at least after replacing $\shL$ by $\shL^{\otimes 2}$.  Over infinite ground fields, such open sets contain
rational points, so Lefschetz pencils indeed exist.  This carries over to finite ground fields $k=\FF_{p^\nu}$, by passing to 
higher tensor powers of $\shL$, as explained in  \cite{Poonen 2009} and \cite{Nguyen 2005}. 

The blowing-up $\tilde{S}=\Bl_Z(S)$ is another smooth surface,
endowed with a Lefschetz fibration  $\tilde{f}:\tilde{S}\ra\PP^1$. The fibers $\tilde{f}^{-1}(t)$ are the strict transforms
of the $C_t$, and the projection $\tilde{f}^{-1}(t)\ra C_t$ are isomorphisms, because the center $Z$ is Cartier on $C_t$.
In particular, no exceptional curve  for the blowing-up $\tilde{S}\ra S$ belongs to a fiber of $\tilde{f}$.

We now form the cartesian square 
$$
\begin{CD}
X_n	@>h_n>>	\tilde{S}	\\
@Vf_nVV		@VV\tilde{f}V \\
\PP^1	@>>F^n>	\PP^1,
\end{CD}
$$
where $F^n$ denotes the iterated relative Frobenius map for the projective line.
Note that the horizontal arrows are  finite universal homeomorphisms that are flat of degree $p^n$,
and that $\Sing(X_n/k)$ maps to the finite set $\Sing(\tilde{S}/\PP^1)$.
By Serre's Criterion,  the surface $X_n$ is normal.

\begin{proposition}
\mylabel{prop:frobenius base-change rdp}
Let $x\in X_n$ be a $k$-rational point mapping to $\Sing(\tilde{S}/\PP^1)$. Then the local ring $\O_{X_n,x}$ is a rational double point
of type $A_l$ with $l=p^n-1$. 
Moreover, the tangent sheaf $\Theta_{X_n/k}$ is locally free,
and the evaluation pairing $\Theta_{X_n/k}\otimes\Omega^1_{X_n/k}\ra\O_{X_n}$ is surjective.
\end{proposition}

\proof
Let $\tilde{s}\in\tilde{S}$ and $t\in\PP^1$  be the images of $x\in X_n$. Let also $R=\O_{\tilde{S},\tilde{s}}^\wedge$, $A=\O_{\PP^1,z}^\wedge$ the completions of the corresponding local rings and  $z\in A$ be a uniformizer. By the definition of Lefschetz fibration, the fiber $\tilde{f}^{-1}(t)$ has a simple normal crossing singularity at $\tilde{s}\in \tilde{S}$. Therefore, $\mathfrak{m}_R/(z)=(\bar{u},\bar{v})$ and $\bar{u}\bar{v}=0$, where $\mathfrak{m}_R$ is the maximal ideal of $R$, and $\bar{u}, \bar{v}$ are the images in $R/(z)$ of two elements $u,v \in \mathfrak{m}_R$. In particular, $uv=zw$, for some $w\in R$.  Considering now that $R$ is a two-dimensional regular local ring, 
\[
\dim \mathfrak{m}_R/ \mathfrak{m}_R^2 =\dim \mathfrak{m}_R/( \mathfrak{m}_R^2+(z))=2.
\]
This implies that $\mathfrak{m}_R=(u,v)$, and $z \in \mathfrak{m}_R^2$. Hence $u,v$ is a set of regular system of parameters for $R$. Therefore, since $uv\in \mathfrak{m}_R^2-\mathfrak{m}_R^3$, it follows that $w\in R^{\ast}$ is a unit. Then by replacing $u$ with $uw$, $R$ is described by the equation $z=uv$. Hence, the complete local ring at $x\in X$  is described by $z^{p^n}=uv$, which yields a rational double point of type $A_{p^n-1}$, as claimed.
From Proposition  \ref{prop:rdp of a-type} below we obtain that  the tangent sheaf $\Theta_{X_n/k}$ is locally free
and the evaluation pairing  $\Theta_{X_n/k}\otimes\Omega^1_{X_n/k}\ra\O_{X_n}$ is surjective.
\qed

\medskip
The group scheme $G_n=\PGL_2[F^n]$ acts on the projective line $\PP^1$ in a canonical way,
and the iterated relative Frobenius map $F^n:\PP^1\ra \PP^1$ obviously factors over the quotient $\PP^1/G_n$. 

\begin{lemma}
\mylabel{lem:quotient on projective line}
The induced map $\PP^1/G_n\ra\PP^1$ is an isomorphism.
\end{lemma}

\proof
By induction, it suffices to treat the case $n=1$. The map in question is a map of degree one between
proper regular curves, hence must be an  isomorphism.
\qed

\begin{lemma}
\mylabel{lem:frobenius base-change general type}
Suppose $S$ is a minimal surface of general type. Then the
fibers of  $\tilde{f}:\tilde{S}\ra\PP^1$ are stable curves of some genus $g\geq 2$, and the base-change $X$
is a RDP surface of general type.
\end{lemma}

\proof
\newcommand{\tD}{\tilde{D}}
We may assume that $k$ is algebraically closed. Recall that the fibers $\tilde{f}^{-1}(t)\subset \tilde{S}$
are copies of the curves $C_t\subset S$. The  irreducible fibers have genus $g\geq 2$,  because $\omega_\tS$ is ample on each movable curve.
Suppose now that   $C_t=D\cup D'$ is reducible, and write $\tD,\tD'\subset\tS $ for the strict transforms of the 
irreducible components. From 
$$
0=(\tD+\tD')^2=\tD^2 + 2(\tD\cdot \tD') + \tD'^2=  \tD^2 + 2 + \tD'^2
$$
we infer $\tD^2=\tD'^2=-1$, and thus $D^2,D'^2\geq -1$.
The Adjunction Formula gives $-2\chi(\O_D)=(\omega_S\cdot D)+D^2\geq 0-1$, and thus $h^1(\O_D)\geq 1/2$,
and likewise for $D'$. It follows that  $C_t$ is a stable curve of genus $g\geq 2$.

Let $\omega_{F^n}$ be the relative dualizing sheaf   for $F^n:\PP^1\ra \PP^1$. Then  
\[
\omega_{F^n}=\omega_{\mathbb{P}^1}\otimes ((F^n)^{\ast}\omega_{\mathbb{P}^1})^{-1}=\omega_{\mathbb{P}^1} \otimes \omega_{\mathbb{P}^1}^{\otimes-p^n}=\mathcal{O}_{\mathbb{P}^1}(-2+2p^n),
\]
according to \cite{Kleiman 1980}, Corollary 24. Consequently 
$\omega_{X_n}= h^*(\omega_{\tilde{S}}) \otimes f_n^*(\O_{\PP^1}(-2+2p^n))$.
The first tensor factor is big, the second effective, and it follows that the tensor product  remains big. 
Using that  the singularities of $X_n$ are rational double points, we infer on the minimal resolution of singularities
$X'_n\ra X_n$, the dualizing sheaf stays big. Summing up, $X_n$ is a RDP surface of general type.
\qed
 
\medskip\noindent
\emph{Proof of Theorem \ref{thm:existence with higher height}.}
By construction, the canonical morphism $X_n\ra S$ is an
alteration whose function field extension is purely inseparable of degree  $p^n$.
From the above observation, we see that in 
$$
X_n=\tilde{S}\times_{\PP^1}(\PP^1,F^n)=\Bl_Z(S)\times_{\PP^1}(\PP^1,F^n),
$$
the   $G_n$-action on the second factor induces an action on the fiber product,
and the projection $f_n:X_n\ra \PP^1$ is equivariant. The action is faithful on the projective line,
so the same holds on $X_n$.  
The statement on the singularities and the tangent sheaf follow from  Proposition \ref{prop:frobenius base-change rdp}.
If $S$ is a minimal surface of general type, then $X_n$ is  a RDP surface of general type, by Lemma \ref{lem:frobenius base-change general type}.
\qed

\medskip 
Let us also make the following observation:

\begin{proposition}
\mylabel{prop:frobenius base-change inertia}
For each   $x\in X_n$, the inertia group scheme $G_x$ inside the base-change of $G=\PGL_2[F]$ has a two-dimensional
Lie algebra.
\end{proposition}

\proof
We may assume that $x$ is a rational point, and it suffices to verify the corresponding statement for 
the image $b=f_n(x)$ on the projective line.
Then the Lie algebra in question is the vector space 
\[
H^0(\PP^1,\Theta_{\PP^1/k}(-b))=H^0(\PP^1,\Theta_{\PP^1/k}(-1))=H^0(\PP^1,\mathcal{O}_{\PP^1}(1)),
\]
which is two-dimensional.
\qed

\medskip
In odd characteristics the canonical projection $\SL_2\ra\PGL_2$ is an isomorphism, By taking this into consideration and also from~\cite[Theorem 12.1]{Schroeer; Tziolas 2023}, the above results can be summarized as follows:

\begin{proposition}
\mylabel{prop:frobenius base-change for p odd}
For $p\neq 2$ and $n\geq 1$, the RDP surfaces $X_n$ of general type carries a fixed point free $\liesl_2$-triple and $\mathrm{Aut}_{X/k}[F] \cong \mathrm{SL}_2[F]$.
\end{proposition}
 

 \section{Rational double points and tangent modules}
\mylabel{sec:rational double points}

We   encountered RDP surfaces  having  a fixed point free $\liesl_2$-triple, where furthermore the tangent sheaf is 
locally free and the evaluation pairing with K\"ahler differentials is surjective. 
In light of the classification of rational double points, it is natural to ask which of them have   one or both of these  properties.
Such questions showed up in many circumstances (for example \cite{Lipman 1965}, \cite{Wahl 1975},  \cite{Ekedahl; Hyland; Shepherd-Barron 2012}, 
\cite{Schroeer 2008}, \cite{Hirokado 2019}, \cite{Schroeer 2021}, \cite{Graf 2022}, \cite{Kawakami 2022}, \cite{Liedtke 2024}, \cite{Martin; Stadlmayr 2024}), 
and the following observation shed more light on these issues.

Let us start with some homological algebra:   Suppose $R$ is a local noetherian ring
with maximal ideal $\maxid$ and residue field $k=R/\maxid$. Recall that each finitely generated module
$M$ has  \emph{free resolutions}
$$
\ldots\lra R^{\oplus r_3}\stackrel{\varphi_3}{\lra} R^{\oplus r_2}\stackrel{\varphi_2}{\lra} 
R^{\oplus r_1}\stackrel{\varphi_1}{\lra} R^{\oplus r_0}\lra M\lra 0.
$$
Such a resolution is called \emph{minimal} if in  $\varphi_i$ viewed as a matrix the entries belong to the maximal ideal, 
or equivalently $\varphi_i(R^{\oplus r_i}) \subset \maxid^{\oplus r_{i-1}}$.
Then the ranks $r_i\geq 0$ depend only on $M$, and are called \emph{Betti numbers} $ b_i(M)=r_i$.
Obviously, the  \emph{ projective dimension} $\pd(M)\leq \infty$ is finite  if and only if one of the Betti numbers vanishes.
In what follows  we write
$$
M^{\ast}=\Hom(M,R)\quadand E=\Ext^1(M,R),
$$
and seek to understand when    $M^{\ast}$ is free,
or the evaluation pairing $M^{\ast}\otimes M\ra R$ given by $\psi\otimes a\mapsto \psi(a)$ is surjective, in dependence
of   Betti numbers of the involved modules.
 
\begin{lemma}
\mylabel{lem:homological characterizations} 
Suppose  $b_2(M)=0$. Then the  following holds:
\begin{enumerate}
\item The  dual  module $M^{\ast}$ is free if and only if $b_3(E)=0$.
\item The evaluation paring  $M^{\ast}\otimes M\ra R$ is surjective  if and only if $b_1(E)<b_0(M)$.
\end{enumerate} 
\end{lemma}

\proof 
Choose a minimal free resolution $0\ra R^{\oplus s}\stackrel{{}^tJ}{\ra}  R^{\oplus r}\ra M\ra 0$.
The resulting  long exact sequence
$$
0\lra M^{\ast} \lra R^{\oplus r }\stackrel{J}{\lra} R^{\oplus s}\lra E\lra 0
$$
already gives (i).  For the remaining statement, note that 
the entries in the matrices ${}^tJ$ and $J$   belong to the maximal ideal, hence
$$
b_0(M)=r \quadand  b_1(M)=s=b_0(E).
$$
For assertion (ii), suppose first that $b_1(E)<r$. Without loss of generality,
we may assume that all but the last standard basis vectors  $e_1,\ldots,e_{r-1}\in R^{\oplus r}$ already generate
the image of $J$. Then $R^{\oplus r-1}\times 0\stackrel{J}{\ra} R^{\oplus s}\ra E\ra 0$
is the beginning of a  free resolution. Write $N$ for the syzygy module, and consider the commutative diagram
$$
\begin{CD}
0	@>>>	N	@>>>	R^{\oplus r-1}	@>>>	\Image(J)	@>>> 0\\
@.		@VVV		@VVV			@VV\id V\\
0	@>>>	M^{\ast}	@>>>	R^{\oplus r}	@>>> 	\Image(J)	@>>> 0.	
\end{CD}
$$
The Snake Lemma provides a surjection $\psi:M^{\ast} \ra R^{\oplus r}/R^{\oplus r-1}=R$ with kernel $N$. For the image  
$a\in M$  of the last basis vector $e_m\in R^{\oplus m}$ we then have $\psi(a)=1$.

Conversely, suppose there are  $\psi\in M^{\ast}$, $a\in M$ such that  $\psi(a)$ is a unit.
Then our original module $M$ splits off an invertible summand, and takes the form  $M=M^{\prime}\oplus R$.
Clearly $M^{\prime}$ and $M$ have the same projective dimension, and the Betti numbers are related by
$b_0(M^{\prime})=b_0(M)-1$ and $ b_1(E)=b_1(E')=b_0(M^{\prime})<b_0(M)$.
\qed

\medskip
There are  further practical characterizations for freeness of $M^{\ast}$ or surjectivity of $M^{\ast}\otimes M\ra R$
in the case where $R$ is a local complete intersection ring, $E$ is a Cohen Macauley $R$-module and 
\begin{equation}
\label{eq:assumptions for numerical characterizations}
b_2(M)=0\quadand b_1(M)=1\quadand \dim \mathrm{Supp}(E)=0 \quadand p>0.
\end{equation} 
Then the minimal resolution takes
the from $0\ra R\stackrel{{}^tJ}{\ra} R^{\oplus r}\ra M\ra 0$. In the ensuing  exact sequence $0\ra M^{\ast}\ra R^{\oplus r}\stackrel{J}{\ra} R\ra E\ra 0$,
the image of $J$ is an $\maxid$-primary ideal $\ideala=(f_1,\ldots,f_r)$, and we write $\ideala^{[p]}=(f_1^p,\ldots,f_r^p)$ for its 
\emph{Frobenius power}. 
Now the three \emph{modules of finite length} $ R/\ideala$ and $R/\ideala^{[p]}$ and $\ideala/\ideala\maxid$  become relevant:
 
\begin{proposition}
\mylabel{prop:numerical characterizations}
Suppose that $R$ is a complete intersection ring and the assumptions \eqref{eq:assumptions for numerical characterizations} hold. Then we have  the following numerical characterizations:
\begin{enumerate}
\item The  dual   module $M^{\ast}$ is free if and only if  $\dim (R)=2$ and $\length(R/\ideala^{[p]})=p^2\cdot  \length(R/\ideala)$.
\item The evaluation  pairing $M^{\ast}\otimes M\ra R$ is surjective if and only if the inequality
$\length(\ideala/\ideala\maxid)< r$ holds.
\end{enumerate} 
Moreover, if the local ring $R$ is two-dimensional,  the surjectivity
of the pairing $M^{\ast}\otimes M\ra R$ implies the freeness of the module $M^{\ast}$.
\end{proposition}

\proof
We start with (i). According to \cite{Miller 2003}, Corollary 5.2.3 we have
$$
\pd(E)<\infty \quad\Longleftrightarrow \quad \length(R/\ideala)=p^d\cdot\length(R/\ideala^{[p]}),
$$
a fact already used in \cite{Schroeer 2008}, Section 4. If the above equivalent conditions hold, the
Auslander--Buchsbaum Formula  (\cite{Eisenbud 1995}, Theorem 19.9) gives
$$
\pd(M^{\ast})=\pd(E)-2=\depth(R)-\depth(E)-2= \depth(R)-2,
$$
by the assumptions $(9)$. So  $M^{\ast}$ is free if and only if  $\depth(R)=2$ and also $\length(R/\ideala)=p^d\cdot\length(R/\ideala^{[p]})$. But since $R$ is a complete intersection ring, it is Cohen Macauley and therefore, since $\depth(R)=2$, $d=\dim (R)=2$.

Suppose the evaluation pairing is surjective. By Lemma \ref{lem:homological characterizations},
the ideal $\ideala=(f_1,\ldots,f_r)$ is already generated by less than $r$ elements,
and the Nakayama Lemma yields $\length(\ideala/\maxid\ideala)<r$.  The converse is likewise.

Finally, assume that $\dim(R)=\depth(R)=2$, and that the equivalent conditions of (ii) hold. Then $r=3$, and the ideal $\ideala$ is already generated by $r-1=2$ element. Write 
$\ideala=(f,g)$. By assumption  $E\cong R/\ideala$ has zero-dimensional support,
Hence the closed sets $V(f ),V(g)\subset\Spec(R)$ have no common irreducible component.
Using that $R$ is Cohen--Macaulay, we infer that $f,g$ form a regular sequence, so the Koszul complex
$$
0\lra \Lambda^2(F)\stackrel{{}^t(g,-f)}{\lra} \Lambda^1(F)\stackrel{(f,g)}{\lra} \Lambda^0(F)\lra R/\ideala\lra 0
$$
for $F=R^{\oplus 2}$ reveals $\pd(R/\ideala)\leq 3$, hence $M^{\ast}$ is free.
\qed

\medskip
The above  observations  from homological algebra apply in the following geometric situation: 
Fix a ground field $k$, for the moment of arbitrary characteristic $p\geq 0$,
and consider   local noetherian rings of the form 
$$
R=S^{-1}k[T_1,\ldots,T_m]  /(P_1,\ldots,P_n),
$$
where the multiplicative system $S$ comprises all polynomials having non-zero constant terms,
and $P_1,\ldots,P_n\not\in S$.
With $\ideala=(P_1,\ldots,P_n)$ and $A=k[T_1,\ldots,T_m]$ we get an exact sequence
\begin{equation}
\label{eq:sequence with kaehler differentials}
\ideala/\ideala^2\lra \Omega^1_{A/k}\otimes_A R\lra \Omega^1_{R/k}\lra 0.
\end{equation} 
where $\ideala/\ideala^2$ is   generated by the classes of $P_i$, and the map on the left is given by $P_i+\ideala^2\mapsto dP_i\otimes 1$.

If  the $P_1,\ldots,P_m$ form a regular  sequence in $S^{-1}k[T_1,\ldots,T_m]$, the module $\ideala/\ideala^2$ is free and the $P_i$
provide a basis.  If furthermore $R$ is geometrically reduced, the map to $\Omega^1_{A/k}\otimes_A R$ injective,
and the exact sequence \eqref{eq:sequence with kaehler differentials} becomes a free resolution
$0\ra R^{\oplus n}\ra R^{\oplus m}\ra \Omega^1_{R/k}\ra 0$. It is minimal provided in addition that $\mathrm{edim}(R)=m$.
Summing up, Lemma \ref{lem:homological characterizations} and Proposition \ref{prop:numerical characterizations} 
apply for $M=\Omega^1_{R/k}$ under the above three assumptions,  and give some information on the   modules  
$$
\Theta_{R/k}=T^0_{R/k}=\Hom(\Omega^1_{R/k},R)\quadand T^1_{R/k}=\Ext^1(\Omega^1_{R/k},R).
$$
Note that the latter can be seen as the space of first order deformations for the scheme $U=\Spec(R)$,
as explained in \cite{Artin 1976}.

Suppose now that the ground field $k$ is algebraically closed, and recall that  the   \emph{rational double points}
have been classified \cite{Artin 1977} and can be described by explicit equations $P(x,y,z)=0$.
The dual graph of the exceptional divisor on the minimal resolution of singularities corresponds 
to the \emph{Dynkin diagrams} $A_l\; (l\geq 1)$ or  $D_l\; (l\geq 4)$   or $E_6,E_7,E_8$. This already determines
the formal isomorphism class, at least in characteristics  $p\geq 7$. 
For the $E$-types at the remaining primes and   the $D$-types in characteristic two
one has to introduce    \emph{upper indices} to distinguish
formal isomorphism class; these reflect the deformability of the singularity. 

\begin{proposition}
\mylabel{prop:rdp of a-type}
For the rational double points of type $A_l$, $l\geq 1$ we have:
$$
\text{$\Theta_{R/k}$ is free}\quad\Longleftrightarrow\quad
\text{$\Theta_{R/k}\otimes\Omega^1_{R/k}\ra R$ is surjective}\quad\Longleftrightarrow\quad
\text{$l\equiv -1$ modulo $p$.}
$$
Moreover, in characteristic  $p\geq 7$ there are no other rational double points where
the tangent module   is free or the paring is surjective,
and for $p\geq 3$ there are at least no such of $D$-type.
\end{proposition}

\proof 
In what follows, $P=P(x,y,z)$ is  the polynomial defining the rational double point, $\idealb$ the ideal in $A=k[x,y,z]$  
generated by $P$ and its partial derivative, and  $\idealb'$ is  generated by $P$ and the $p$-th powers
of the partial derivatives.

For type $A_l$  the four polynomials in question  are 
$$
P=z^{l+1}-xy,\quad P_z= (l+1)z^l,\quad  P_x= -y,\quad  P_y= -x.
$$
Suppose first $l+1\equiv 0$ modulo $p$.
Then $\ideala=\idealb R=(x,y)$, and from Lemma \ref{lem:homological characterizations} we see that the pairing $\Theta_{R/k}\otimes \Omega^1_{R/k} \ra R$ is surjective,
while  Proposition \ref{prop:numerical characterizations} ensures that $\Theta_{R/k}$ is free.
Suppose now that $l+1\not\equiv 0$. Then $\length(A/\idealb)=l$.   On the other hand,  
$\bar{A}=k[x,y,z]/(P,x^p,y^p)$  is free as module over $k[x,y]/(x^p,y^p)$, and $1,z,\ldots,z^l$ provide a basis,
giving $\length(\bar{A})=p^2\cdot l$. In light of Proposition \ref{prop:numerical characterizations}, it remains
to check that $z^{lp}\in\bar{A}$ is non-zero. Indeed, otherwise $l+1$ divides $lp$, which implies $l+1=p$, contradiction.

Suppose next that $p\neq 2$. For   type $D_l$, $l\geq 4$ we have  
$$
P=z^2+x^2y+xy^{l-1},\quad P_z= 2z,\quad  P_x= 2xy+y^{l-1},\quad  P_y= x^2+(l-1)xy^{l-2}.
$$
Inside the polynomial ring $A=k[x,y,z]$, we consider the ideals $\idealb=(P,P_z,P_x,P_y)$ and $\idealb'= \idealb^{[p]}+(P)$.
Our task is to show $\length(A/\idealb')\neq p^2\cdot\length(A/\idealb)$.
First observe that $(l-1)xP_x-yP_y=2(l-1) (1-y)x^2$. It follows that $(P_z,P_x,P_y)$ is also generated by
$$
 Q_1=z\quadand Q_2=x^2 ,\quadand Q_3= 2xy+y^{l-1},
$$
and thus $\length(A/\idealb)=2(l-1)$. Likewise we see that $\bar{A}=A/(P,Q_2^p,Q_3^p)$ has length $p^2\cdot 2(l-1)$.
It thus suffices to verify that the class of $Q_1^p=z^p$ in $\bar{A}$ does not vanish.
The latter is freely generated, as module over $k[x]/(x^{2p})$, by the monomials $y^iz^j$ with 
$0\leq i< (l-1)p$ and $0\leq j<2$. Computing modulo $P$, we see
$$
z^p = (z^2)^{(p-1)/2} \equiv (-x^2y-xy^{l-1})^{(p-1)/2} = \pm x^{(p-1)/2}y^{(l-1)(p-1)/2} + (...),
$$
where the remaining summand have lower degree in $y$. It follows that $z^p\in \bar{A}$ is indeed non-zero.

It remains to treat type $E_6$, $E_7$,  $E_8$ in characteristic $p\geq 7$, where the  respective polynomials  are
$P=z^2+x^3+y^4$ and $ P=z^2+x^3+xy^3$ and $P=z^2+x^3+y^5$. The arguments are as above, and left to the reader.
\qed

\medskip
As usual, the small primes are in need of  special attention.
There are only finitely many rational double points of $E$-type; using computer algebra \cite{Magma},
we compute the   lengths for the relevant modules $ R/\ideala, R/\ideala^{[p]}$
together with  $\ideala/\ideala\maxid$ occurring  in Proposition \ref{prop:numerical characterizations}. 
The results are collected in table \ref{tab:e in small char}, and immediately yield the following:

\begin{proposition}
\mylabel{prop:rdp of e-type}
The rational double points of $E$-type  where the tangent module $\Theta_{R/k}$ is free
are precisely for  
$$
E_8^0\;(p=5)\quadand E_6^0,E_7^0,E_8^0\; (p=3)\quadand E_6^0,E_7^0,E_7^1,E_7^2, E_8^0,E_8^1, E_8^2\; (p=2).
$$
The evaluation pairing $\Theta_{R/k}\otimes\Omega^1_{R/k}\ra R$ is surjective precisely for   
$E_8^0\;(p=5)$ and $  E_6^0,E_7^0,E_8^0\; (p\leq 3)$.
\end{proposition}

The hardest challenge is   to understand the rational double points of  $D$-type   in characteristic two.
Using computer algebra  to compute  examples, one immediately comes up with the following:

\begin{proposition}
\mylabel{prop:rdp of d-type}
In characteristic $p=2$, the rational double points of type $D_l^r$, $l\geq 4$  have the following behavior with regards to the tangent module  $\Theta_{R/k}$ and
the  pairing $\Theta_{R/k}\otimes\Omega^1_{R/k}\ra R$:
\begin{enumerate}
\item For $r=0$, the  tangent module is free and the pairing is surjective.  
\item For $l=2n$ even and $1\leq r\leq n-1$,  the tangent module is free but the pairing is not surjective.
\item For $l=2n+1$ odd and $1\leq r\leq n-1$, neither the tangent module is  free nor  the pairing is   surjective.
\end{enumerate}
\end{proposition}

\proof
For $D_{2n}^0$ and $D_{2n+1}^0$ the respective polynomials are
$P=z^2+x^2y+xy^n$ and $P=z^2+x^2y+y^nz$. One of the partial derivatives $P_z$ or $P_x$ vanishes,
so (i) holds by Lemma \ref{lem:homological characterizations}.
 
Suppose now $r\geq 1$ and $l=2n$. Then the partial derivatives of the  defining polynomial $P=z^2 +x^2y+xy^n+xy^{n-r}z$
admit factorizations 
$$
P_z=x\cdot y^{n-r}\quadand  P_x=y^{n-r}\cdot(y^r+ z)\quadand P_y = x\cdot Q 
$$
with $Q=x+ny^{n-1}+(n-r)y^{n-r-1}z$,  and thus can be seen as maximal minors of the $3\times 2$-matrix
$$
\varphi_1={}^t\begin{pmatrix}
y^r+z\	& x	& 0\\
Q	& 0	& y^{n-r}
\end{pmatrix}.
$$
Extending the matrix with $\varphi_0=(P_z,P_x,P_y)$ as first row and applying Laplace Expansion, we get a complex
$0\ra R^{\oplus 2}\stackrel{\varphi_1}{\ra} R^{\oplus 3}\stackrel{\varphi_0}{\ra} \ideala\ra 0$
for the ideal $\ideala\subset R$   generated by the partial derivatives. According to the Hilbert--Burch Theorem
(\cite{Eisenbud 1995}, Theorem 20.15), this complex is actually exact, hence $\Theta_{R/k}$ is free.
All entries of $\varphi_1$ belong to the maximal ideal, so Lemma \ref{lem:homological characterizations} tell us that 
$\Theta_{R/k}\otimes\Omega^1_{R/k}\ra R$ is not surjective. This settles (ii). 

It remains to treat the case $r\geq 1$ and $l=2n+1$. Now the defining polynomial is 
$P=z^2+x^2y+y^nz+xy^{n-r}z$, with partial derivatives
$$
P_z=y^n+xy^{n-r},\quad  P_x=y^{n-r}z,\quad  P_y=x^2+ny^{n-1}z+ (n-r)xy^{n-r-1}z.
$$
To simplify notation set $\nu=n-r$, and consider the following three derivations
\begin{equation}
\label{eq:three derivations}
\begin{gathered}
\delta_1=
z\frac{\partial}{\partial z} + (x+y^r)\frac{\partial}{\partial x},\\
\delta_2=
(\nu+1)y^{n-r}z\frac{\partial}{\partial z} + (\nu y^n+z)\frac{\partial}{\partial x} + y^{n-r+1}\frac{\partial}{\partial y},\\
\delta_3= 
(xz+y^rz+\nu y^{n-r-1}z)\frac{\partial}{\partial z} +  (y^{2r}+(\nu+1)y^{n-1}z)\frac{\partial}{\partial x} +  y^{n-r}z\frac{\partial}{\partial y}\\
\end{gathered}
\end{equation}
of $A=k[x,y,z]$. One easily checks that each of them stabilizes the principal ideal generated by $P$,
and thus give rise to elements $\bar{\delta}_1,\bar{\delta}_2,\bar{\delta}_3\in\Theta_{R/k}$. 
Let $J_2\in\Mat_3(A)$ be the coefficient matrix for the above derivations, 
and $J_1=(P_z,P_x,P_z)$  be the Jacobi matrix.  Then
$$
J_1\cdot J_2\equiv 0 \mod P,
\qquad
J_2\equiv\begin{pmatrix}
z	& 0	& 0\\
x	& z	& 0\\
0	& 0	& 0\\
\end{pmatrix}
\mod \maxid_A.
$$
Seeking a contradiction, we now assume that $\Theta_{R/k}$ is free.
Suppose for the moment that $b_1(T^1_{R/k})=3$. Then $R^{\oplus 3}\stackrel{J_1}{\ra}R\ra T^1_{R/k}\ra 0 $ is the beginning of a minimal resolution,
and for all $\bar{\delta}$ from $\Theta_{R/k}\subset R^{\oplus 3}$ the coordinates belong to $\maxid_R$.
Using that the first two rows in $J_2$ are linearly independent modulo $\maxid_A^2$,
we infer that $\bar{\delta}_1,\bar{\delta}_2\in\Theta_{R/k}$ form a basis. We thus have
$\bar{\delta}_3=Q_1\bar{\delta}_1 +Q_2\bar{\delta}_2$ for some polynomials $Q_i=Q_i(x,y,z)$.  
Comparing coefficients in the third coordinate, we arrive at  
$y^{n-r}z \in (y^{n-r+1},P)$.
On the other hand, the residue class ring $k[x,y,z]/(y^{n-r+1},P)$ is   free as $k[x]$-module, 
with basis given by monomials $y^iz^j$ with $0\leq i< n-r+1$ and $0\leq j<2$.
So $y^{n-r}z$ is a basis member,  contradiction.

It remains to verify   $b_1(T^1_{R/k})=3$. Seeking a contradiction, we assume that the ideal $\idealb=(P,P_z,P_x,P_y)$
in $A=k[x,y,z]$ is   generated by $P$ together with only two  partial derivatives.
To use  Gr\"obner bases techniques, introduced weights 
$$
\deg(x)=\deg(z)=n-1\quadand \deg(y)=1,
$$
and define a \emph{monomial ordering}  $x^iy^rz^j\succ x^{i'}y^{r'}z^{j'}$
if $\deg(x^iy^rz^j)< \deg(x^{i'}y^{r'}z^{j'})$, and in case of equality use the lexicographic ordering with $x\succ z\succ y$.
By design, the  leading terms are    as follows:
\begin{equation}
\label{eq:leading terms}
P=z^2+(\ldots),\quad P_z=y^n+(\ldots),  \quad P_x=y^{n-r}z, \quad P_y=x^2+(\ldots).
\end{equation} 
We claim that the above generators form a \emph{standard basis}.
Recall that this notion was introduced, in one form or another,  by Hironaka
\cite{Hironaka 1964}, Buchberger \cite{Buchberger 1970}  and Grauert \cite{Grauert 1972}, and is now usually seen 
as a    variant of \emph{Gr\"obner bases} in the setting
of localized or completed polynomial rings. See for example \cite{Greuel; Pfister 2002} or \cite{Cox; Little; O'Shea 2005} for more details.
To apply the  \emph{Buchberger Criterion} we have to compute
$S$-polynomials and check that division with reminder reduces them to zero: This indeed holds for
$$
S(P_z,P_x)=xP_x\quadand S(P,P_x)=(x^2y+y^n+xy^{n-r})P_x.
$$
The other $S$-polynomials also reduce to zero because  the leading terms of the generators are pairwise prime,
(\cite{Cox; Little; O'Shea 2007}, Section 2, \S 9, Proposition 4).  Thus \eqref{eq:leading terms} from a standard basis, hence
$$
\length(A/\idealb)=\length(k[x,y,z]/(z^2,y^n,y^{n-r}z,x^2))= 4n-2r,
$$
as stated without proof in \cite{Artin 1977}. Writing $\idealb_z=(P,P_x,P_y)$ and  $\idealb_x=(P,P_z,P_y)$ and $\idealb_y=(P,P_z,P_x)$, 
we likewise see that the generators form a standard basis, and obtain
$$
\length(A/\idealb_z)=\length(A/\idealb_y)=\infty \quadand \length(A/\idealb_x)=4n.
$$
In turn, one cannot remove in  $\idealb=(P,P_z,P_x,P_y)$ a partial derivative to obtain a smaller generating set,
and thus $b_1(T^1_{R/k})=3$.
\qed

\begin{table} 
{\footnotesize
$$
\begin{array}{*{12}{l}}
\toprule
	& E_6^0		& E_6^1		& E_7^0		& E_7^1		& E_7^2		& E_7^3		& E_8^0		& E_8^1		& E_8^2		& E_8^3		& E_8^4\\		
\toprule
p=2	& 8,32  	& 6,28		& 14,56		& 12,48		& 10,40		& 8,35		& 16,64		& 14,56		& 12,48		& 10,44		& 8,37 \\		
 	& 2  		& 3		& 2		& 3		& 3		& 3		& 2		& 3		& 3		& 3 		& 3 \\		

\midrule
p=3	& 9,81 		& 7,71		& 9,81		& 7,75		& 		& 		& 12,108	& 10,99		& 8,85\\
	& 2 		& 3		& 2		& 3		& 		& 		& 2		& 3		& 3\\

\midrule
p=5	& 6, 173	&		& 7, 198 	&		&		&		& 10, 250 	& 8, 239\\
	& 3 		& 		& 3		& 		& 		& 		& 2		& 3		& 		&  		&  \\		
\bottomrule
\end{array}
$$
}
\medskip
\caption{Lenghts of $R/\ideala, R/\ideala^{[p]}$ and $\ideala/\ideala\maxid$ for RDP of $E$-type  for  $p\leq 5$.}
\label{tab:e in small char}
\end{table} 

\begin{remark}
Note that the statements pertaining to the freeness of $\Theta_{R/k}$ already appears in the work of Graf~\cite{Graf 2022}, with an argument relying on a computer program. This is fine in each individual case, but we could not see how this settles infinite families like $D_n$. Also note that Matsumoto (\cite{Matsumoto 2023a}, Proposition 4.7 and \cite{Matsumoto 2023b}, Theorem 3.3) obtained similar results on the pairing $\Theta_{R/k}\otimes\Omega^1_{R/k}\ra R$. The difference is that he assumed that surjectivity is obtained with p-closed rather than an arbitrary derivations and he did not relate this with the freeness of the tangent sheaf. Our arguments are more general and they can be applied to treat other classes of singularities too.
\end{remark}


\section{Rational double points with \texorpdfstring{$\liesl_2$}{sl2}-triples}
\mylabel{sec:rational double points with triples}

Let $k$ be an algebraically closed ground field of characteristic $p>0$.
In this final section we take up the natural question
which rational double points admit $\liesl_2$-triples. This is  of local nature, and we work in the following setting:
Let $Y$ be a   connected normal surface that is separated and of finite type, endowed with  an $\liesl_2$-triple,
in other words, a non-trivial action of the height-one group scheme $G=\SL_2[F]$.
Suppose $y\in Y$ is a RDP singularity and that there are no further singularities, and write  $f:X\ra Y$ for the minimal resolution of singularities.
In characteristic two the following arguments are somewhat    inconclusive, and for  the sake  of exposition  we assume  $p\neq 2$ throughout.

\begin{proposition}
\mylabel{prop:possible rdp with triple}
The rational double point $y\in Y$ is of type $A_1$, or belongs to the following list:
$$
A_l\;(l\equiv -1)\quadand E_8^0\;(p=5)\quadand E_6^0, E_7^0,E_8^0,\; (p=3).
$$
Here  the above congruence is modulo $p$. Moreover,  the $G$-action on $Y$ does not extend to $X$, except for type $A_1$.
\end{proposition}

\proof
Let $E_1,\ldots,E_l\subset X$ be the exceptional divisors.
According to \cite{Artin 1966}, Theorem 4 the schematic fiber $Z=f^{-1}(y)$ coincides with the fundamental cycle of the singularity.
For type $A_l$ the cycle is reduced, while  for  $D_l$ and $E_6,E_7,E_8$ the multiplicities in $Z=\sum_{i=1}^l m_iE_i$  are as follows:
\begin{gather*}
\dynkin[labels={1,2,2}]A3 \cdots \dynkin[labels={ ,2,1,1}]D4\quadand
\dynkin[labels={1,2,2,3,2,1}]E6\qquad\dynkin[labels={2,2,3,4,3,2,1}]E7\qquad\dynkin[labels={2,3,4,6,5,4,3,2}]E8
\end{gather*}
 For $l=1$ our singularity is of type $A_1$, and there is nothing to prove. Suppose now $l\geq 2$. 
Since $p\neq 2$,   we find by inspection in each case  a pair  of indices $r\neq s$ such that  $E_r\cap E_s$ 
is non-empty and  $p\nmid m_rm_s$.

Seeking a contradiction, we suppose that the $G$-action on $Y$ lifts to $X$. Then the fiber $Z=f^{-1}(y)$ is $G$-stable (\cite{Martin 2022}, Corollary 2.18).
In turn, the  $m_iE_i$  are $G$-stable (\cite{Brion; Schroeer 2023}, Lemma 2.3), and the same holds for $D=m_rE_r\cup m_sE_s$. 
It then follows that $D_\red=E_r\cup E_s$ is $G$-stable 
(\cite{Tziolas 2022}, Proposition 4.2), therefore the unique point $x\in E_r\cap E_s$ is $G$-fixed.
If the $G$-action on $E_r=\PP^1$ is non-trivial, the induced map $G\ra\PGL_2[F]$ is an isomorphism, and we see that $x$ is not fixed,
contradiction. So $E_r$ is $G$-fixed, and  the same holds for $E_s$. In particular, the fixed scheme $X^G$ contains
the singular curve $E_r\cup E_s$.  On the other hand, for the multiplicative  group schemes $\mu_p\subset G$ 
the fixed scheme $X^{\mu_p}$ is smooth (\cite{Abramovich; Temkin 2018}, Proposition 5.1.16), contradiction. 

This establishes that the $G$-action on $Y$ does not lift to the minimal resolution of singularities. Consequently,  the canonical map $f_*(\Theta_{X/k})\ra\Theta_{Y/k}$ is not bijective.
According to Hirokado's result (\cite{Hirokado 2019}, Theorem 1.1), this holds precisely when the singularity $y\in Y$ is of the following types. $A_l\;(l\equiv -1), \;E_8^0\;(p=5),\; E_6^0, E_6^1, E_7^0,E_8^0, E_8^1\; (p=3).$ In order to conclude the proof of the proposition, it suffices to show that the cases $E_6^1$ and $E_8^1$ in characteristic $p=3$ do not happen.

Suppose that $y\in Y$ is a singularity of type $E_6^1$ or $E_8^1$. Let $g\colon Y^{\prime}\rightarrow  Y$ be the blow up of $X$ at $y$ and $E\cong \mathbb{P}^1_k$ the $g$-exceptional curve. A straightforward local calculation shows that $Y^{\prime}$ has a single singular point $Q$. Moreover, if $y\in Y$ is of type $E_6^1$, $Q\in Y^{\prime}$  is a rational double point of type $A_5$. If $y\in Y$ is of type $E_8^1$, then $Q\in Y^{\prime}$ is of type $E_7^0$. According to Proposition~\ref{prop:rdp of e-type}, the evaluation pairing in both cases is not surjective and therefore the singular point $y\in Y$ is a fixed point of every vector field of $Y$. Therefore, the natural map $g_{\ast}\Theta_{Y^{\prime}}\rightarrow \Theta_Y$ is an isomorphism.

Let now $\mathfrak{sl}_2(k)\subset H^0(\Theta_Y)$ be the embedding of $\mathfrak{sl}_2(k)$ in $H^0(\Theta_Y)$ corresponding to the $\mathfrak{sl}_2$-triple. Since $g_{\ast}\Theta_{Y^{\prime}}=\Theta_Y$, this lifts in $H^0(\Theta_{Y^{\prime}})$. Since the $g$-exceptional curve $E$ is $D$-invariant for every vector field $D$ of $Y^{\prime}$, there exists a map of restricted Lie algebras $\tau\colon H^0(\Theta_{Y^{\prime}})\rightarrow H^0(T_E)$.  We will show that the image of $\tau$ is one-dimensional, in particular $\tau$ is not surjective. This implies that the restriction of $\tau$ on $\mathfrak{sl}_2(k)$, gives a map of restricted Lie algebras $\mathfrak{sl}_2(k) \rightarrow H^0(T_E)$ which is not injective, since $\mathfrak{sl}_2(k)$ is $3$-dimensional. Considering that $\mathfrak{sl}_2(k)$ is a simple Lie algebra, this means that this map is actually zero. Hence the lifting of the $\mathfrak{sl}_2$-triple in $Y^{\prime}$ vanishes on $E$. Hence it fixes the singular point $Q\in Y^{\prime}$ and therefore lifts to the blow up 
$Y^{\prime\prime}\rightarrow Y^{\prime}$ of $Q\in Y^{\prime}$. A straightforward calculation shows that $Y^{\prime\prime}$ has a single singular point $R$ over $Q\in Y^{\prime}$. If $Q\in Y^{\prime}$ is of type $A_5$, then $R\in Y^{\prime\prime}$ is of type $A_4$ and if $Q\in Y^{\prime}$ is of type $E_7^0$, then $Q\in Y^{\prime}$  is of type $D_5$. But now in both cases, according to~\cite[Theorem 1.1]{Hirokado 2019}, every vector field of $Q\in Y^{\prime\prime}$ lifts to its  minimal resolution. Therefore, the original $\mathfrak{sl}_2$-triple of $Y$ lifts to the minimal resolution $X$ of $Y$. But this is impossible from the previous arguments.

It remains to show that the map $\tau\colon H^0(\Theta_{Y^{\prime}})\rightarrow H^0(T_E)$ has one-dimensional image. This is local at  the singular point $y\in Y$. We will do only the $E_6^1$ case. The $E_8^1$ is identical and is omitted. 

Suppose that $y\in Y$ is of type $E_6^1$. Then
\[
R=\hat{\mathcal{O}}_{Y,y}\cong \frac{k[[x,y,z]]}{(x^2+y^3+z^4+y^2z^2)}.
\]
 By using the computer algebra program Singular~\cite{DGPS} we have found that $\Theta_{R/k}$ is generated by the derivations 
 \[
 \delta_i=f_{1,i} \frac{\partial}{\partial x}+f_{2,i}\frac{\partial}{\partial y}+f_{3,i}\frac{\partial}{\partial z},
\]
for $i=1,2,3,4$, where
\[
(f_{i,j})=
\begin{pmatrix}
0                      &    xz            &  yz^2   &y^2z-z^3   \\
y^2-z^2            &z-yz          &  -x      & 0              \\
-yz                   &  y -z^2       & 0          & -x           \\
\end{pmatrix}
\]
Straightforward explicit local calculations in the blow up $Y^{\prime}\rightarrow Y$ using the local generators $\delta_i$, $i=1,\ldots, 6$ of $\Theta_{R/k}$ above, show that the image of the map $\tau\colon H^0(\Theta_{Y^{\prime}})\rightarrow H^0(T_E)$ is one-dimensional.
\qed

\begin{remark}
The argument regarding the exclusion of the $E_6^1$ and $E_8^1$ cases was communicated to us by Yuya Matsumoto. We are grateful to him for his help.
\end{remark}
\medskip
The rings  $A=k[u,v_1,\ldots,v_n]/\ideala $ where the ideal is generated by polynomials
of the form $P=u^n-f(v_1,\ldots,v_n)$ with $n\equiv 0$  admit  an obvious $\liesl_2$-triple, given by
the derivations
$h=2u\partial/\partial u$  and $e=u^2\partial/\partial u$ and $ f=\partial/\partial u$, a fact already used in Section \ref{sec:lefschetz pencils}.
In light of Artin's classification \cite{Artin 1977} of the rational double points in terms of equations $P(x,y,z)=0$, 
we immediately see that the types 
\begin{equation}
\label{eq:obvious types}
A_l\;(l\equiv -1)\quadand E_8^0\;(p=5)\quadand E_6^0, E_7^0,E_8^0\; (p=3) 
\end{equation} 
admit $\liesl_2$-triples. For   $A_1$ the equation is $z^2-xy=0$, and then  
$$ 
h=-2x\frac{\partial}{\partial x} + 2y\frac{\partial}{\partial y}\quadand
e=2z\frac{\partial}{\partial x} + y\frac{\partial}{\partial z}  \quadand
f=2z\frac{\partial}{\partial y} + x\frac{\partial}{\partial z}
$$
yields an  $\liesl_2$-triple. 

Write $R=\O_{Y,y}$ for the local ring of the rational double point $y\in Y$.
If the $G$-action is fixed-point free, the evaluation pairing $\Theta_{R/k}\otimes\Omega^1_{R/k}\ra R$ is surjective.
Conversely, suppose that the evaluation pairing is surjective.
According to Proposition \ref{prop:rdp of a-type} and Proposition \ref{prop:rdp of e-type}, this happens precisely for the types
listed in \eqref{eq:obvious types}.
A priori,    $y\in Y$   might be  $G$-fixed. Strangely, this is only possible in two cases:

\begin{corollary}
\mylabel{cor:pairing surjective point fixed}
It the evaluation pairing  $\Theta_{R/k}\otimes\Omega^1_{R/k}\ra R$ is surjective and the singularity $y\in Y$ is   $G$-fixed,
the  type is   $E_6^0,E_8^0\; (p=3)$ .
\end{corollary}

\proof
The $G$-action lifts to the blowing-up $Y'=\Bl_y(Y)$. This is a normal surface whose singularities are rational double points.
It turns out that there is at most one singularity $y'\in Y'$, and the types are related as follows:
$$
\begin{array}{lllll}
\toprule
Y	& A_l	& E_6 	& E_7	& E_8\\
\midrule
Y'	& A_{l-2}	& A_5	& D_6	& E_7\\ 
\bottomrule
\end{array}
$$
The minimal resolution $X\ra Y$ factors over $Y'$, and by our Proposition, the $G$-action  does not lift to $X$.
Applying Hirokado's result (\cite{Hirokado 2019}, Theorem 1.1) to the   morphism $X\ra Y'$, we see that we must have $p=3$
and $y'\in Y'$ has type $A_5$ or $E_7^0$. Our assertion follows.
\qed


\end{document}